\documentclass[11pt]{amsart}
\usepackage{amsxtra}
\usepackage{amssymb}
\usepackage{mathtools}
\theoremstyle{plain}
\newcommand{\cleqn}{\setcounter{equation}{0}}
\newcommand{\clth}{\setcounter{theorem}{0}}
\newcommand {\sectionnew}[1]{\section{#1}\cleqn\clth}
\newcommand{\nn}{\hfill\nonumber}
\newtheorem{theorem}{Theorem}[section]
\newtheorem{lemma}[theorem]{Lemma}
\newtheorem{definition-theorem}[theorem]{Definition-Theorem}
\newtheorem{proposition}[theorem]{Proposition}
\newtheorem{corollary}[theorem]{Corollary}
\newtheorem{definition}[theorem]{Definition}
\newtheorem{example}[theorem]{Example}
\newtheorem{remark}[theorem]{Remark}
\newtheorem{conjecture}[theorem]{Conjecture}

\newcommand \bth[1] { \begin{theorem}\label{t#1} }
\newcommand \ble[1] { \begin{lemma}\label{l#1} }

\newcommand \bpr[1] { \begin{proposition}\label{p#1} }
\newcommand \bco[1] { \begin{corollary}\label{c#1} }
\newcommand \bde[1] { \begin{definition}\label{d#1}\rm }
\newcommand \bex[1] { \begin{example}\label{e#1}\rm }
\newcommand \bre[1] { \begin{remark}\label{r#1}\rm }
\newcommand \bcj[1] { \begin{conjecture}\label{j#1}\rm }

\renewcommand {\eth} { \end{theorem} }
\newcommand {\ele} { \end{lemma} }

\newcommand {\epr} { \end{proposition} }
\newcommand {\eco} { \end{corollary} }
\newcommand {\ede} { \end{definition} }
\newcommand {\eex} { \end{example} }
\newcommand {\ere} { \end{remark} }
\newcommand {\ecj} { \end{conjecture} }

\newcommand {\enota} { \end{notation} }
\newcommand \thref[1]{Theorem \ref{t#1}}
\newcommand \leref[1]{Lemma \ref{l#1}}
\newcommand \prref[1]{Proposition \ref{p#1}}

\newcommand \deref[1]{Definition \ref{d#1}}

\newcommand \lb[1]{\label{#1}}
\def \Cset {{\mathbb C}}
\def \Aset {{\mathbb A}}
\def \KK {{\mathbb K}}
\def \Zset {{\mathbb Z}}
\def \Nset {{\mathbb N}}
\def \Qset {{\mathbb Q}}

\def \AA  {{\mathcal{A}}}           %mathcal
\def \B  {{\mathcal{B}}}
\def \CC {\mathcal{C}}
\def \FF {{\mathcal{F}}}

\def \NN {{\mathcal{N}}}

\def \VV {{\mathcal{V}}}

\def \LL {{\mathcal{L}}}
 
\def \ZZ {{\mathcal{Z}}}

\def \ZZ {{\mathcal{Z}}}

\def \De {\Delta}   % Greek letters
\def \de {\delta}
\def \al {\alpha}
\def \be {\beta}

\def \La {\Lambda}

\def \ga {\gamma}
\def \de {\delta}

\def \sig {\sigma}

\def \ep {\epsilon}
\def \sig{\sigma}
\def \mt  {\mapsto}
\def \hra {\hookrightarrow}

\def \rcor {\rangle}
\def \lcor {\langle}
\def \ol {\overline}
\def \wt {\widetilde}

\def \id { {\mathrm{id}} }

\DeclareMathOperator \Span { {\mathrm{Span}} }
\DeclareMathOperator \Aut { {\mathrm{Aut}} }

\DeclareMathOperator \charr { {\mathrm{char}} }
\DeclareMathOperator \opp { {\mathrm{op}} }
\DeclareMathOperator \lcm  { {\mathrm{lcm}}  }

\DeclareMathOperator \ad { {\mathrm{ad}} }
\DeclareMathOperator \diag { {\mathrm{diag}} }
\DeclareMathOperator \Ker { {\mathrm{Ker}} }
\DeclareMathOperator \tr { {\mathrm{tr}} }

\DeclareMathOperator \gr  { {\mathrm{gr}} }

\renewcommand \max { {\mathrm{max}} }
\newcommand \Spec { {\mathrm{Spec}} }
\begin{document}
%%%%%%%%%%%%%%%%%%%%%%%%%%%%%%%%%%%%%%%%%%%%%%%%%%%%%%%%%%%%%%%%%%%%%%%%%%%
%%%%%%%%%%%%%%%%%%%%%%    Title    %%%%%%%%%%%%%%%%%%%%%%%%%%%%%%%%%%%%%%%%
\title[Quantized Weyl algebras at roots of unity]
{Quantized Weyl algebras at roots of unity}
\author[Jesse Levitt]{Jesse Levitt}
\address{
Department of Mathematics \\
Louisiana State University \\
Baton Rouge, LA 70803,
USA
}
\email{jlevit3@lsu.edu}
\author[Milen Yakimov]{Milen Yakimov}
\email{yakimov@math.lsu.edu}
\thanks{The research of J.L. has been supported by a GAANN fellowship and a VIGRE fellowship through the NSF grant DMS-0739382 
and that of M.Y. by the NSF grants DMS-1303038 and DMS-1601862.}
%\date{}
\begin{abstract} 
We classify the centers of the quantized Weyl algebras that are polynomial identity algebras and derive explicit formulas 
for the discriminants of these algebras over a general class of polynomial central subalgebras. Two different approaches 
to these formulas are given: one based on Poisson geometry and deformation theory, 
and the other using techniques from quantum cluster algebras.
Furthermore, we classify the PI quantized Weyl algebras that are free over 
their centers and prove that their discriminants are locally dominating and effective. This is applied to solve the 
automorphism and isomorphism problems for this family of algebras and their tensor products.
\end{abstract}
\maketitle
%%%%%%%%%%%%%%%%%%%%   Introduction   %%%%%%%%%%%%%%%%%%%%%%%%%%%%%%%%%%%%%%%%
\sectionnew{Introduction}
\lb{Intro}
\subsection{Quantized Weyl algebras}
\label{1.1}
%##
The quantized Weyl algebras and their generalizations 
have been studied from many different points of view: quantum groups and Hecke type quantizations \cite{GZ,Ma}, structure of prime spectra and representations \cite{Be,G,GL,J}, automorphism and isomorphism problems \cite{BJ,Ga,GH,KL,RS,R}, homological and ring theoretic dimensions \cite{FKK}, quantizations of 
multiplicative hypertoric varieties \cite{Coo,Gan} and others.
Most of these results concern the generic case when the algebras are not polynomial identity (PI).

In this paper we investigate in detail the properties of quantized Weyl algebras as modules over their centers 
in the case when they are PI. 
We classify their centers and then give a classification of the PI quantized Weyl algebras 
that are free over their centers. For those algebras we compute explicitly their discriminants, and prove
that they are locally dominating and effective. 
More generally, we prove a formula for the discriminant of any PI quantized Weyl algebra 
over a central polynomial subalgebra generated by powers of the standard generators of the Weyl 
algebra. Two different proofs are given for the discriminant formulas: one 
based on deformation theory and Poisson geometry, and a second proof using techniques from 
quantum cluster algebras. This is applied to solve the automorphism and isomorphism problems for the family 
of quantized Weyl algebras that are free over their centers, as well as for the larger class of all tensor products 
of such algebras. 

Let $T$ be a integral domain and $n \in \Zset_+$. (Here and below by an integral domain we 
mean a commutative unital ring without zero-divisors.) 
As usual $T^\times$ will denote the group of units of $T$.
Let $E:=(\ep_1, \ldots, \ep_n) \in (T^\times)^n$ and let $B:=(\be_{jk}) \in M_n(T^\times)$ 
be a multiplicatively skew-symmetric matrix; that is $\be_{jk} \be_{kj} =1$ for all $j \neq k$ and 
$\be_{jj}=1$ for all $j \in [1,n]$. (Here and below we set $[1,n] :=\{1, \ldots, n\}$.) The multiparameter 
quantized Weyl algebra 
$A^{E,B}_n (T)$ is the unital associative $T$-algebra with generators
\[
x_1, y_1, x_2, y_2, \ldots, x_n, y_n
\]
and relations
\begin{align}
y_j y_k &= \be_{jk} y_k y_j, \quad & \forall j, k, 
\nn
\\
x_j x_k &= \ep_j \be_{jk} x_k x_j, \quad &j< k,
\nn
\\
x_j y_k &= \be_{kj} y_k x_j, \quad &j<k,
\label{rel}
\\
x_j y_k &= \ep_k \be_{kj} y_k x_j, \quad &j>k,
\nn
\\
x_j y_j - \ep_j y_j x_j &= 1 + \sum_{i=1}^{j-1} (\ep_i -1) y_i x_i, \quad &\forall j.
\nn
\end{align}
For each sequence of positive integers $\chi_1 < \chi_2 < \cdots < \chi_n$, the algebra  
$A^{E, B}_n (T)$ has an $\Nset$-filtration defined by
\[
\deg x_j =  \deg y_j = \chi_j.
\]
The associated graded algebra $\gr A^{E, B}_n (T)$ is the connected $\Nset$-graded  skew polynomial
algebra with generators $\ol{x}_1, \ol{y}_1, \ldots, \ol{x}_n, \ol{y}_n$ (sitting in degrees $\chi_1, \chi_1, \ldots, \chi_n, \chi_n$)
and relations as above with the exception of the last one which is replaced by 
\[
\ol{x}_j \ol{y}_j = \ep_j \ol{y}_j \ol{x}_j,  \quad \forall j.
\]

The algebra $\gr A^{E, B}_n (T)$ is PI if and only if all $\ep_j, \be_{jk} \in T^\times$ are roots of unity, in which case $\gr A^{E, B}_n (T)$
is module finite over its center. In this situation $A^{E, B}_n (T)$ is module finite over its center if and only if in addition $\ep_j \neq 1$ for all $j$
in the case when $\charr T =0$. 

{\em{Throughout the paper we will assume that}}
\begin{equation}
\label{assum}
\mbox{$\ep_j, \be_{jk} \in T^\times$ are roots of unity and $\ep_j \neq 1$ for all $j$.}
\end{equation}
%%%%%%%%%%%%%%%
\subsection{Results}
\label{1.2}
Consider a family $\FF$ of filtered PI algebras $A$ whose associated graded algebras $\gr A$ are skew polynomial algebras.
Chan, Young and Zhang \cite{CYZ}, based on \cite{CPWZ1,BZ}, proposed the following program for studying the 
family $\FF$:
\begin{enumerate}
\item Classify the algebras $A$ in the family $\FF$ for which $\gr A$ is free over its center $\ZZ(\gr A)$. For those algebras:
\item Describe the center $\ZZ(A)$ in terms of $\ZZ(\gr A)$.
\item Obtain an explicit formula for the discriminant $d(A/\ZZ(A))$.
\item Classify $\Aut(A)$ and obtain other related properties of $A$ (if applicable), e.g. Zariski cancellative, etc.
\end{enumerate}

Here and below, for an algebra $B$, we denote by $\ZZ(B)$ its center.

We carry out this program for the family of quantized Weyl algebras, assuming \eqref{assum}. Let
\begin{equation}
\label{scalars}
\ep_j = \exp( 2 \pi \sqrt{-1} m_j/d_j), \; \; \be_{jk} = \exp( 2 \pi \sqrt{-1} m_{jk}/d_{jk})
\end{equation}
for some $m_j, d_j, d_{jk} \in \Zset_+$ and $m_{jk} \in \Zset$
such that $\gcd (m_j, d_j) =1$, $\gcd (m_{jk}, d_{jk}) =1$ and $d_{jk} = d_{kj}$, $ m_{jk} = - m_{kj}$.
We allow $\be_{jk}=1$ (that is $m_{jk}=0$) for any choice of $j,k$. Denote
\[
D(E,B) = \lcm(d_j, d_{jk}, 1 \leq j,k \leq n).
\]
The imaginary exponents in \eqref{scalars} are interpreted as follows. We assume that $T$ contains a $D(E,B)$-th primitive root of unity, 
to be denoted by $\exp( 2 \pi \sqrt{-1}/D(E,B))$. In particular $\charr T \nmid D(E,B)$.
The imaginary exponents in \eqref{scalars} are the corresponding powers of this element 
of $T^\times$.

The first problem in part (1) of the program was solved 
by Chan, Young and Zhang for arbitrary skew polynomial algebras in terms of a deep number theoretic 
condition on the entries of the defining matrix for the algebra \cite[Theorem 0.3]{CYZ}.
This condition is elaborate, and even checking if a given square integer matrix satisfies 
it, is nontrivial. 

For the quantized Weyl algebras we prove a more explicit criterion that solves the first part of the program and, furthermore, we prove 
that $\ZZ(A^{E,B}_n(T)) \cong \ZZ(\gr A^{E,B}_n(T))$:
\medskip
\\
\noindent
{\bf{Theorem A.}} {\em{Let $T$ be an integral domain. Assume the setting of \eqref{scalars}. 
\begin{enumerate}
\item[(i)] The canonical map $\gr \colon A^{E,B}_n(T) \to \gr A^{E,B}_n(T)$ induces a $T$-algebra isomorphism
$\ZZ(A^{E,B}_n(T)) \cong \ZZ(\gr A^{E,B}_n(T))$ if $\ep_j-1 \in T^\times$ for all $j$. 

\item[(ii)] 
For all integral domains $T$, the following are equivalent for the quantized Weyl algebra $A^{E,B}_n(T)$:
\begin{enumerate}
\item The algebra $\gr A^{E,B}_n(T)$ is free over its center;
\item The center $\ZZ(\gr A^{E,B}_n(T))$ is a polynomial algebra;
\item The algebra $A^{E,B}_n(T)$ is free over its center;
\item The center $\ZZ(A^{E,B}_n(T))$ is a polynomial algebra;
\item $d_j | d_l$ and $d_{jk} |d _l$ for all $j \leq l$ and $k \in [1,n]$.
\end{enumerate}
If these conditions are satisfied, then 
\begin{equation}
\label{center}
\ZZ(A^{E,B}_n(T)) = T[x_j^{d_j}, y_j^{d_j}, 1 \leq j \leq n].
\end{equation}
\end{enumerate}
}}
In addition, in \thref{Centers} we give a description of $\ZZ(A^{E,B}_n(T))$ for all PI quantized Weyl algebras $A^{E,B}_n(T)$.
We note that the conditions in Theorem A (ii) are satisfied in the important uniparameter 
case when $\ep_1 = \ldots =\ep_n$ and $\be_{jk}=1$ for all $j,k$.

Next, we investigate the algebras $A^{E,B}_n(T)$ satisfying the conditions in Theorem A (ii), solving parts (2) and (3) of the above program. 
We offer two very different approaches based on 
\begin{enumerate}
\item[(A1)] deformation theory and Poisson geometry, and 
\item[(A2)] quantum cluster algebras, 
\end{enumerate}
respectively. When the conditions in Theorem A (ii) are satisfied, $D(E,B)$ is given by
\begin{equation}
\label{dd}
D(E,B) = d_n. 
\end{equation}
Let
\begin{equation}
\label{XYgen}
X_j:= x_j^{d_j}, Y_j := y_j^{d_j} \in \ZZ(A^{E,B}_n(T)).
\end{equation}
Define recursively the sequence of elements $Z_0, \ldots, Z_n \in \ZZ(A^{E,B}_n(T))$ by $Z_0 := 1$
and
\begin{equation}
\label{Zgen}
Z_j := - (1 - \ep_j)^{d_j} Y_j X_j +Z_{j-1}^{d_j/d_{j-1}} \quad \mbox{for} \; \; j \in [1,n].
\end{equation}
In \prref{Poisson} we show that $Z_j = z_j^{d_j}$ in terms of the frequently used normal elements
\[
z_j := 1 + (\ep_1 -1) y_1 x_1 + \cdots + (\ep_j -1) y_j x_j = [x_j,y_j] \in A^{E, B}_n (T).
\]

The algebras $A^{E,B}_n(T)$ are specializations of algebras over $T[q^{\pm1}]$. This induces a canonical Poisson algebra structure on 
$\ZZ(A^{E,B}_n(T))$ and, more generally, turns $A^{E,B}_n(T)$ into a Poisson order over its center \cite{BrGo}.
Our next result describes the induced Poisson structure on $\ZZ(A^{E,B}_n(T))$ and the set of Poisson prime elements of this algebra
(see \S 2.1-2.2 for background). We compute the discriminant of $A^{E,B}_n(T)$ over $\ZZ(A^{E,B}_n(T))$ 
and establish dominating/effectiveness properties for it (see \S 2.4 for definitions). 
\medskip
\\
 \noindent
{\bf{Theorem B.}} {\em{Let $T$ be an integral domain. In the setting of \eqref{scalars}, assume that the conditions in 
Theorem A (ii) are satisfied and that $\charr T \nmid d_n = D(E,B)$.
\begin{enumerate}
\item[(i)] The induced Poisson structure on $\ZZ(A^{E,B}_n(T))$ is given by \eqref{bracxy}. 
If $T$ is a field of characteristic 0, 
then the only Poisson prime elements of $\ZZ(A^{E,B}_n(T))$ {\em{(}}up to associates{\em{)}} are $Z_1, \ldots, Z_n$.

\item[(ii)] The discriminant of $A^{E,B}_n(T)$ over its center is given by
\begin{align*}
d(A^{E,B}_n(T)/ \ZZ(A^{E,B}_n(T))) &=_{T^\times} \eta Z_1^{N^2 (d_1-1)/d_1} \ldots Z_n^{N^2 (d_n-1)/d_n}
\\
&=_{T^\times}  \eta z_1^{N^2 (d_1-1)} \ldots z_n^{N^2 (d_n-1)},
\end{align*}
where $N:=d_1 \ldots d_n$ and 
$\eta := \Big( N \prod_{j=1}^n ([d_j-1]_{\ep_j}!) \Big)^{N^2} \in T$, see \eqref{eta2} for an alternate expression for $\eta$. 

\item[(iii)] Let $A^{E_1,B_1}_{n_1}(T), \ldots, A^{E_l,B_l}_{n_l}(T)$ be a collection of quantized Weyl algebras 
satisfying the conditions in Theorem A (ii) and $A$ be their tensor product over $T$. 
If 
\[
\charr T \nmid D(E_1, B_1) \ldots D(E_l, B_l), 
\]
then
the discriminant $d(A, \ZZ(A))$ is locally dominating
and effective in the sense of \cite{CPWZ1} and \cite{BZ}.
\end{enumerate}
}}

The case of the theorem when $n=n_1 = \ldots = n_l=1$ was proved in \cite{CYZ}.

As usual, we set $[k]_q := (1-q^k)/(1-q)$ and $[k]_q! = [1]_q \ldots [k]_q$.
In part (ii) the discriminant is computed using the trace $\tr \colon A^{E,B}_n(T) \to \ZZ(A^{E,B}_n(T))$ associated 
to the embeddings $A^{E,B}_n(T) \hra M_{N^2} (\ZZ(A^{E,B}_n(T)))$ from $\ZZ(A^{E,B}_n(T))$-bases of $A^{E,B}_n(T)$.
Here and below for two elements $r_1$ and $r_2$ of an integral domain $R$ we write $r _1 =_{R^\times} r_2$ 
if $r_1$ and $r_2$ are associates, i.e., if $r_1 = r_2 u$ for some unit $u$ of $R$. 

Furthermore, in \thref{general-discr} we prove a general formula for the discriminant
\[
d(A^{E,B}_n(T), T[x_j^{L_j}, y_j^{L_j}, 1 \leq j \leq n])
\]
for any collection of central elements $x_j^{L_j}, y_j^{L_j} \in \ZZ(A^{E,B}_n(T))$ without assuming that the conditions in 
Theorem A (ii) are satisfied. We do not state the formula here because it is considerably more complicated than the one in Theorem B (ii). The proof also works for 
the more general case when the powers of $x_j$ and $y_j$ are different but the discriminant formula is even more involved. This is why we restrict
ourselves to equal powers of $x_j$ and $y_j$.

We give two proofs of Theorem B (ii) based on the two above mentioned approaches.
Both proofs not only yield the explicit formulas for the discriminants, 
but also explain the intrinsic meaning of the factors of the discriminants from two different perspectives. In fact, the 
two proofs concern the irreducible factors in two different algebras: the first one deals with the factors in the 
centers of the PI quantized Weyl algebras and the second with the factors in the whole PI quantized Weyl algebras:
\begin{enumerate}
\item The first proof relies on the geometry of the induced Poisson structure from the specialization. 
The factors of the discriminants are precisely {\em{the unique Poisson primes}} $Z_1, \ldots Z_n$ from Theorem B (i).
This proof 
works for base rings $T$ of arbitrary characteristics (using certain base change and filtration arguments), even though some parts 
of the Poisson geometric setting require $T=\Cset$.
\item The second proof uses quantum cluster algebra structures on the quantized Weyl algebras. It shows that the irreducible 
factors of these discriminants in the whole quantized Weyl algebras $A^{E,B}_n(T)$ are precisely {\em{the frozen variables of these quantum
cluster algebras}}.
\end{enumerate}

We combine Theorems A and B and the results in \cite{CPWZ1,BZ} to resolve the automorphism and isomorphism 
problems for the tensor products of the family of quantized Weyl algebras $A^{E,B}_n(T)$ satisfying the conditions 
in Theorem A (ii). 
\medskip
\\
\noindent
{\bf{Theorem C.}}  {\em{Let $A = A^{E_1,B_1}_{n_1}(T) \otimes_T \cdots \otimes_T A^{E_l, B_l}_{n_l}(T)$ 
for a collection of quantized Weyl algebras 
over an integral domain $T$, satisfying the conditions in Theorem A (ii). Assume that 
$\charr T \nmid D(E_1, B_1) \ldots D(E_l, B_l)$ and recall \eqref{dd}.

{\em{(i)}} If $\phi \in \Aut_T(A)$, then the following hold:
\begin{enumerate}
\item[(1)] There exists $\sig \in S_l$ such that $n_i = n_{\sig(i)}$ and $\phi(A^{E_i,B_i}_{n_i}(T)) = A^{E_{\sig(i)}, B_{\sig(i)}}_{n_{\sig(i)}}(T)$
for all $i \in [1,l]$.
\item[(2)] For a given $i \in [1,l]$, denote the standard generators of $A^{E_i,B_i}_{n_i}(T)$ and $A^{E_{\sig(i)}, B_{\sig(i)}}_{n_{\sig(i)}}(T)$
by $x_1, y_1, \ldots, x_n, y_n$ and $x'_1, y'_1, \ldots, x'_n, y'_n$ where $n = n_i = n_{\sig(i)}$. There exist scalars
$\mu_1, \nu_1, \ldots, \mu_n, \nu_n \in T^\times$ and a sequence $(\tau_1, \ldots, \tau_n) \in \{\pm 1 \}^n$ such that
\begin{align*}
&\phi(x_j)= \mu_j x'_j, \quad \phi(y_j)= \nu_j y'_j, &\mbox{if} \; \; \tau_j =1, \\
&\phi(x_j)= \mu_j y'_j, \quad \phi(y_j)= \nu_j x'_j, &\mbox{if} \; \; \tau_j=-1.
\end{align*}
The scalars satisfy the following equalities for $B_i=(\be_{jk})$, $B_{\sig(i)} = (\be'_{jk})$, $E_i = (\ep_1, \ldots, \ep_n)$ and
$E_{\sig(i)} = (\ep'_1, \ldots, \ep'_n)$:
\begin{align}
& \mu_j \nu_j = \tau_j \prod_{
\begin{matrix}
1 \leq k \leq j \\
\tau_k = -1
\end{matrix}} \ep_k^{-1}, \quad \ep'_j = \ep_j^{\tau_j}, \quad \forall j,
\label{ident1}
\\
&\be'_{jk} = 
\begin{cases}
\be_{jk}^{\tau_j}, &\mbox{if} \; \; \tau_k =1
\\
(\ep_j \be_{jk})^{- \tau_j}, & \mbox{if} \; \; \tau_k =-1,
\end{cases}
\quad\quad \forall j<k.
\label{ident2}
\end{align}
\end{enumerate}

{\em{(ii)}} Every map $\phi$ on the $x$- and $y$-generators of $A$ with the above properties extends to a $T$-linear automorphism 
of $A$.

{\em{(iii)}} The algebra $A$ is strongly Zariski cancellative and ${\mathrm{LND}}^H$-rigid, see \S \ref{2.4} for terminology.}}
\medskip

The case of the theorem when $n_1 = \ldots = n_l=1$ was proved in \cite{CYZ}. 

The theorem also solves the isomorphism problem for the family of algebras consisting of all tensor products of collections of quantized Weyl algebras 
$A^{E,B}_n(T)$ satisfying the conditions in Theorem A (ii). Indeed, if $A$ and $B$ are two such tensor products and $\psi \colon A \stackrel{\cong}\to B$ is 
an algebra isomorphism, then $\phi := \psi \otimes \psi^{-1} \in \Aut_T (A \otimes_T B)$ and $A \otimes_T B$ is one of the algebras whose 
automorphisms are classified by the theorem. This way one recovers all isomorphisms $\psi \colon A \stackrel{\cong}\to B$.

The theorem can be specialized to classify the automorphisms of a single quantized Weyl algebra $A^{E,B}_n(T)$ and the isomorphisms
between two given quantized Weyl algebras $\psi \colon A^{E_1,B_1}_{n_1}(T) \stackrel{\cong} \to A^{E_2,B_2}_{n_2}(T)$ satisfying the conditions in Theorem A (ii). 
We state those results in Sect. \ref{5.4}.
\medskip
\\
\noindent
{\bf Acknowledgements.} We are grateful to Ken Goodearl for his very helpful comments and for suggesting an improvement 
of our original proof of Theorem B (i) that extended the result from the field of complex numbers to all fields of characteristic 0. 
We would also like to thank Bach Nguyen for his very valuable suggestions. 
%%%%%%%%%%%%%%%%%%%%%%%%%%%%%%%%%%%%%%
\sectionnew{Poisson geometry of specializations and discriminants}
\lb{Discr}
This section contains background material on Poisson structures for algebras 
obtained by specialization, discriminants of noncommutative algebras, their 
relations to Poisson geometry and their applications to the automorphism and isomorphism 
problems for algebras.
\subsection{Poisson structures for specializations}
\label{2.1}
Let $T$ be an integral domain and $A$ be a $T[q^{\pm1}]$-algebra which is a torsion free $T[q^{\pm 1}]$-module.
Given $\ep \in T^\times$, the specialization of $A$ at $\ep$ is
the $T$-algebra $A_\ep: = A/(q -\ep)A$. The natural projection $\sig \colon A \to A_\ep$ is a homomorphism of $T$-algebras.
The specialization gives rise to a canonical Poisson structure on $\ZZ(A_\ep)$ and a lifting of the hamiltonian derivations of $\ZZ(A_\ep)$ 
to derivations of $A_\ep$:

(i) The canonical structure of Poisson algebra on $\ZZ(A_\ep)$ is defined as follows. For $z_1, z_2 \in \ZZ(A_\ep)$, choose $c_i \in \sig^{-1}(z_i)$ and set 
\[
\{ z_1, z_2 \} := \sig \big( (c_1 c_2 -c_2 c_1)/(q-\ep) \big).
\]
One easily verifies that the RHS does not depend on the choice of preimages. We have $\{z_1, z_2 \} \in \ZZ(A_\ep)$ because
\[
[\{ z_1, z_2 \}, \sig (a)] = \sig ( [[c_1,c_2],a]/(q-\ep))=0, \quad \forall a \in A.
\]

(ii) Given $z \in \ZZ(A_\ep)$, one constructs \cite{DKP1,Ha} lifts of the hamiltonian derivation $x \mt \{ z, x \}$ of $( \ZZ(A_\ep), \{. , .\})$ to derivations 
of $A_\ep$ as follows. Choose $c \in \sig^{-1}(z)$ and set
\[
\partial_c ( \sig(a)) := \sig \big( [c,a]/(q-\ep) \big).
\]
The fraction $[c,a]/(q-\ep)$ is an element of $A$ because $[\sig(c), \sig(a)]=0$, so $[c,a] \in \ker \sig = (q-\ep) A$. One easily verifies that the 
RHS does not depend on the choice of preimage of $\sig(a)$ and that $\partial_c$ is a derivation of $A_\ep$. The derivations corresponding 
to two different preimages of $z$ differ by an inner derivation of $A_\ep$:
\[
\partial_c- \partial_{c'} = \ad \sig((c-c')/(q-\ep)), \quad \forall c,c' \in \sig^{-1}(z).
\]

The above situation was axiomatized by Brown and Gordon \cite{BrGo} as the notion of Poisson order. In this language, the constructions 
(i) and (ii) say that $A_\ep$ is a Poisson $C_\ep$-order for each Poisson subalgebra $C_\ep$ of $\ZZ(A_\ep)$ 
such that $A_\ep$ is a finite rank $C_\ep$-module.
%%%%%%%%%%
\subsection{Poisson prime elements}
\label{2.2}
Assume that $(P, \{.,.\})$ is a Poisson algebra which is an integral domain as an algebra.

\bde{prime-norm} (i) An element $a \in P$ is called {\em{Poisson normal}} if there exists a 
Poisson derivation $\partial$ of $P$ such that 
\[
\{a, x \} = a \partial(x) \quad \forall x \in P.
\]
Equivalently, $a \in P$ is Poisson normal if and only if the ideal $(a)$ is Poisson.

(ii) An element $p \in P$ is called {\em{Poisson prime}} if it is a prime element of the commutative algebra $P$ 
which is Poisson normal. Equivalently, $p \in P$ is Poisson prime if and only if $(p)$ is a nonzero prime and Poisson
ideal of $P$.
\ede
If $P$ is the coordinate ring of a smooth complex affine Poisson variety $W$, then $f \in \Cset[W]$ is Poisson 
prime if and only if $f$ is prime and its zero locus $\VV(f)$ is a union of symplectic leaves of $W$, see \cite[Remark 2.4 (iii)]{NTY}.
%%%%%%%%%%
\subsection{Discriminants and their relation to Poisson algebras}
\label{2.3}
Let $A$ be an associative algebra and $C$ be a subalgebra of $\ZZ(A)$. A $C$-valued trace on $A$ is 
a $C$-linear map $\tr \colon A \to C$ such that 
\[
\tr(x y) = \tr (yx) \quad \forall x,y \in A.
\]
If $A$ is a free $C$-module of finite rank, one defines \cite{Re,CPWZ1} the discriminant $d(A/C) \in C$ as follows. For every two $C$-bases 
$\B$ and $\B'$ of $A$ set
\[
d_N(\B, \B' \colon \tr) := \det \big( [\tr(b b')]_{b \in \B, b' \in \B'} \big),
\]
where $N = |\B|=|\B'|$.
If $\B_1$ and $\B'_1$ are two other $C$-bases of $A$, then \cite[Eq. (1.10.1)]{CPWZ2}
\[
d_N(\B_1, \B'_1 \colon \tr) =_{C^\times} d_N(\B, \B' \colon \tr).
\]
We define the discriminant of $A$ over $C$ by 
\[
d(A/C):=_{C^\times} d_N(\B, \B' : \tr)
\]
for any $C$-bases $\B$ and $\B'$ of $A$. More generally, if $A$ is a finite rank $C$-module which is not necessarily free, there are 
notions of {\em{discriminant}} and {\em{modified discriminant ideals}} of $A$ over $C$, \cite{Re,CPWZ1,CPWZ2}. They are the ideals of $C$ generated by the 
elements of the form $d_N(\B,\B : \tr)$ and $d_N(\B, \B' : \tr)$, respectively. We also set
\[
d_N(\B : \tr):= d_N(\B, \B : \tr).
\]

When $A$ is free of finite rank over the subalgebra $C \subset \ZZ(A)$, there is a natural $C$-valued trace map on $A$, often 
called the {\em{internal}} trace of $A$. Any $C$-basis of $A$ gives rise to an embedding  $A \hra M_N(C)$ where 
$N$ is the rank of $A$ over $C$. The internal trace of $A$ is the composition of this embedding with 
the standard trace $M_N(C) \to C$.
\bth{prod} \cite[Theorem 3.2]{NTY} Let $A$ be a $\KK[q^{\pm1}]$-algebra for a field $\KK$ of characteristic 0 
which is a torsion free $\KK[q^{\pm 1}]$-module and $\ep \in \KK^\times$. 
Assume that the specialization $A_\ep := A/(q-\ep)A$ is a free module of finite rank over a Poisson 
subalgebra $C_\ep$ of its center and that $C_\ep$ is a unique factorization domain as a commutative algebra. 
\begin{enumerate}
\item[(i)] Let $\tr \colon A_\ep \to C_\ep$ be a trace map which commutes with all derivations $\partial$ of $A_\ep$ such that 
$\partial(C_\ep) \subseteq C_\ep$.  The corresponding discriminant
$d(A_\ep/C_\ep)$ either equals $0$ or 
\[
d(A_\ep/C_\ep) =_{C_\ep^\times} \prod_{i=1}^m p_i
\]
for some (not necessarily distinct) Poisson prime elements $p_1, \ldots, p_m \in C_\ep$.
\item[(ii)] The internal trace coming from the freeness of $A_\ep$ as a $C_\ep$-module commutes with all derivations 
$\partial$ of $A_\ep$ such that $\partial(C_\ep) \subseteq C_\ep$. 
\end{enumerate}
\eth
%%%%%%%%%%
\subsection{Applications of discriminants}
\label{2.4}
Discriminants play a major role in many settings in algebraic number theory, algebraic geometry and 
combinatorics \cite{GKZ,S}. In noncommutative algebra they play a key role in the study of orders \cite{Re}.
Recently, many other applications of them were found in the study of the automorphism and isomorphism problems for PI algebras 
and the Zariski cancellation problem for noncommutative algebras \cite{BZ,CPWZ1,CPWZ2}. 

In this subsection we review some terminology and results of Ceken--Palmieri--Wang--Zhang \cite{CPWZ1}, Makar-Limanov \cite{M},
and Bell--Zhang \cite{BZ}.
Fix a unital $T$-algebra $A$. Consider a generating set $x_1, \ldots, x_n$ of $A$, such that $\{1, x_1, \ldots, x_n\}$ 
is $T$-linearly independent. Denote by $\FF_i A$ the filtration of $A$ obtained by assigning 
$\deg x_i =1$, $\forall i$, i.e., $\FF_k (A) := (T. 1 + T x_1 + \cdots + T x_n)^k$. In this setting, we define:
\bde{specialprop} (i) \cite[Definition 2.1 (1)]{CPWZ1}
An element $f \in A$ is called {\em{locally dominating}} if for every $\phi \in \Aut_T(A)$:
\begin{enumerate}
\item $\deg \phi(f) \geq \deg f$ and
\item $\deg \phi(f) > \deg f$ if $\deg(\phi(x_i)) > \deg x_i$ for at least one $i$,
\end{enumerate}
where the degrees are computed with respect to the filtration $\FF_k A$. 

(ii) \cite[Definition 5.1]{BZ} An element $f \in A$ is called {\em{effective}} if $A$ has (a possibly different) filtration $F_k A$ whose associated graded 
algebra is a domain with the following property.
For every testing $\Nset$-filtered PI $T$-algebra $S$, whose associated graded is a domain, 
and for every testing subset $\{y_1, \ldots, y_n \} \subset S$, satisfying
\begin{enumerate}
\item $\{1, y_1, \ldots, y_n\}$ is $T$-linearly independent and
\item $\deg_S y_j \geq \deg_A x_j$ for all $j$ and $\deg_S y_i > \deg_A x_i$ for some $i$ (with respect to the filtration 
$F_k A$),
\end{enumerate}
$f$ has a lift $f(x_1, \ldots, x_n)$ in the free algebra $T \lcor x_1, \ldots x_n\rcor$ 
such that either $f(y_1, \ldots, y_n)=0$ or $\deg_S f(y_1, \ldots, y_n) > \deg_A f$.
\ede
A stronger notion of dominating elements of algebras was defined in \cite{CPWZ1}.

\bde{can} (i) An algebra $A$ is called {\em{cancellative}}, if $A[t] \cong B[t]$ for an algebra $B$ implies $A \cong B$.

(ii) An algebra $A$ is called {\em{strongly cancellative}} if, for all $k \geq 1$, $A[t_1, \ldots, t_k] \cong B[t_1, \ldots, t_k]$
for an algebra $B$ implies $A \cong B$. 
\ede

Denote by ${\mathrm{LND}}(A)$ the $T$-module of locally nilpotent derivations of $A$. The {\em{Makar-Limanov invariant}} \cite{M}
of $A$ is defined by
\[
{\mathrm{ML}}(A) := \bigcap_{\partial \in {\mathrm{LND}}(A)} \Ker \partial,
\]
where one sets ${\mathrm{ML}}(A) := A$ if ${\mathrm{LND}}(A) = 0$.
An algebra $A$ is called \cite{BZ,M} {\em{${\mathrm{LND}}$-rigid}} if ${\mathrm{ML}}(A) = A$, and {\em{strongly ${\mathrm{LND}}$-rigid}} if 
${\mathrm{ML}}(A[t_1, \ldots, t_k]) = A$ for all $k \in \Zset_+$.

A {\em{locally nilpotent higher derivation}} of $A$ is a sequence $\partial :=(\partial_0 = \id, \partial_1, \ldots)$ of $T$-endomorphisms 
of $A$ such that 
\[
a \mt \sum_{j=0}^\infty \partial_j(a) t^k 
\]
is a well defined $T[t]$-algebra automorphism of $A[t]$ (in particular, for all $a \in A$, $\partial_j(a)=0$ for sufficiently large $j$).  
In this case the map $\partial_1$ is necessarily a derivation of $A$.
The set of those elements $\partial$ is denoted by ${\mathrm{LND}}^H(A)$. The {\em{higher Makar-Limanov invariant}} \cite{BZ} of $A$ is 
defined by 
\[
{\mathrm{ML}}^H(A) := \bigcap_{\partial \in {\mathrm{LND}}^H(A)} \Ker \partial, \quad
\mbox{where} \quad \Ker \partial = \bigcap_{j\geq1} \Ker \partial_j. 
\]
The algebra $A$ is called {\em{strongly ${\mathrm{LND}}^H$-rigid}} if 
${\mathrm{ML}}^H(A[t_1, \ldots, t_k]) = A$ for all $k \in \Zset_+$. If $T$ is an extension of $\Qset$,
then $A$ is ${\mathrm{LND}}$-rigid if and only if it is ${\mathrm{LND}}^H$-rigid, \cite[Remark 2.4 (a)]{BZ}.

\bth{appl-dom} Assume that $A$ is a $T$-algebra which is a free module over $\ZZ(A)$ of finite rank.

 {\em{(i) \cite[Theorem 2.7]{CPWZ1}}} If the discriminant $d(A/\ZZ(A))$ is locally dominating with respect 
 to the filtration $\FF_k A$ associated to a set of generators $\{x_1, \ldots, x_n\}$, 
 then every $\phi \in \Aut_T(A)$ is affine in the sense that $\phi(\FF_1 A) = \FF_1 A$. 

{\em{(ii)}} \cite[Theorem 5.2]{BZ} If $A$ is a domain and the discriminant $d(A, \ZZ(A))$ is effective, then
$A$ is strongly ${\mathrm{LND}}^H$-rigid. If, in addition, $A$ has finite GK-dimension, then $A$ is strongly cancellative.
\eth
A stronger cancellation property than the one in (ii) was proved in \cite[Theorem 5.2]{BZ}.
%%%%%%%%%%%%%%%%%%%%%%%%%%%%%%%%%%%%%%
\sectionnew{The centers of the quantized Weyl algebras and their Poisson structures}
\lb{cent-qWeyl}
In this section we describe the centers of the PI quantized Weyl algebras. We classify the PI quantized Weyl algebras
that are free over their centers and relate this property to other properties of the algebras and 
their associated graded algebras. We also compute the induced Poisson structures on their centers when the quantized 
Weyl algebras are realized as specializations of families of algebras.
\subsection{The centers of $A^{E,B}_n(T)$ and $\gr A^{E,B}_n(T)$}
\label{3.1} Throughout this section we assume \eqref{scalars}.
Fix positive integers $\chi_1 < \chi_2 < \cdots < \chi_n$, and consider the $\Nset$-filtration 
on the algebra $A^{E, B}_n (T)$ from \S \ref{1.1} defined by $\deg x_j =  \deg y_j = \chi_j$.
Its filtered components will be denoted by $F_j A^{E, B}_n (T)$, $j \in \Nset$. 
The canonical map $\gr \colon A^{E, B}_n (T) \to \gr A^{E, B}_n (T)$ is given by 
\[
r \mt r + F_{j-1} A^{E, B}_n (T) \quad \mbox{for} \; \;  r \in F_j A^{E, B}_n (T), \; r \notin F_{j-1} A^{E, B}_n (T).
\]
In the notation of \S \ref{1.1}, we have
\[
\ol{x}_j = \gr(x_j), \; \ol{y}_j = \gr(y_j).
\]

As mentioned in \S \ref{1.2} the elements
\[
z_j := 1 + (\ep_1 -1) y_1 x_1 + \cdots + (\ep_j -1) y_j x_j = [x_j, y_j] \in A^{E, B}_n (T), \quad j \in [1,n]
\]
are normal and, more precisely, 
\begin{equation}
\label{z-normal}
z_j x_k = \ep_k^{-\de_{k \leq j}} x_k z_j , \quad 
z_j y_k = \ep_k^{\de_{k \leq j}} y_k z_j, \quad j, k \in [1,n].
\end{equation}
Set $z_0 :=1$. The last relation in \eqref{rel} is transformed into
\begin{equation}
\label{xyz}
x_j y_j = \ep_j y_j x_j + z_{j-1},\; \;  z_{j-1} x_j = x_j z_{j-1}, \;  y_j z_{j-1} = z_{j-1} y_j, \quad j \in [1,n].
\end{equation}
This implies
that
\begin{equation}
\label{xy-normal}
x_j^{d_j} y_j = y_j x_j^{d_j}, \; \; y_j^{d_j} x_j = x_j y_j^{d_j}.
\end{equation}
Thus, $x_j^{d_j}, y_j^{d_j} \in A^{E,B}_n(T)$ normalize all generators $x_k, y_k$. Denote 
\begin{multline*}
C(E,B):=
\{ (b_1, a_1, \ldots, b_n, a_n) \in \Nset^{2n} 
\mid \; \; d_j | (b_j -a_j) , \forall  j \in [1,n], 
\\
\sum_j \frac{(b_j - a_j)m_{jk}}{d_{jk}} + (a_k + \cdots + a_n) \frac{m_k}{d_k} \in \Zset, \forall k \in [1,n] \}.
\end{multline*}
The following proposition is an explicit version of Theorem A (i). Denote
\[
T' := T [(\ep_j -1)^{-1}, 1 \leq j \leq n ].
\]

\bth{Centers} Let $T$ be an integral domain. Assume the setting of \eqref{scalars}. The centers of 
$\gr A^{E, B}_n (T)$ and $A^{E, B}_n (T')$ are given by
\begin{align}
\ZZ( \gr A^{E, B}_n (T)) &= \Span_T \Big\{ \prod \ol{y}_j^{b_j-a_j} (\ol{y}_j \ol{x}_j)^{a_j} \Big\} =
\Span_T \Big\{ \prod \ol{y}_j^{b_j} \ol{x}_j^{a_j} \Big\},
\label{cent-gr}
\\
\ZZ(A^{E, B}_n (T')) &= \Span_{T'} \Big\{ \prod y_j^{\max\{b_j-a_j, 0\}} z_j^{\min\{b_j, a_j\}} x_j^{\max\{a_j-b_j, 0\}} \Big\}, 
\label{cent}
\end{align}
where the spans range over $(b_1, a_1, \ldots, b_n, a_n) \in C(E,B)$.

The $T'$-algebras $\ZZ(A^{E, B}_n (T'))$ and $\ZZ( \gr  A^{E, B}_n (T'))$ are isomorphic.
\eth
The theorem implies that for all integral domains $T$, the center of $A^{E,B}_n(T)$ is given by
\[
\ZZ(A^{E, B}_n (T)) = \Span_{T'} \Big\{ \prod y_j^{\max\{b_j-a_j, 0\}} z_j^{\min\{b_j, a_j\}} x_j^{\max\{a_j-b_j, 0\}} \Big\} 
\bigcap A^{E,B}_n(T).
\]
An important property of the spanning set in this formula is that each monomial involves only 2 factors and 
those factors are powers of the normal elements $y_j^{d_j}$, $z_j$, $x_j^{d_j}$ of $A^{E,B}_n(T)$.
\begin{proof} Let
\[
r:= \ol{y}_1^{b_1} \ol{x}_1^{a_1} \ldots \ol{y}_n^{b_n} \ol{x}_n^{a_n} \in \ZZ( \gr A^{E, B}_n (T)).
\]
The normality \eqref{z-normal} of the elements $z_j$ gives
\[
(\gr z_j) (\gr z_{j-1})^{-1} r = \ep_j^{b_j-a_j} r (\gr z_j) (\gr z_{j-1})^{-1},
\]
thus $d_j | (b_j - a_j)$. The identity $r \ol{y}_k = \ol{y}_k r$ forces
\[
\sum_j \frac{(b_j -a_j) m_{jk}}{d_{jk}} + (a_k + \cdots + a_n) \frac{m_k}{d_k} \in \Zset.
\]
This proves the inclusion $\subseteq$ in \eqref{cent-gr}. The reverse inclusion follows from the facts that
$\gr(z_j) = (\ep_j -1) \ol{y}_j \ol{x}_j$ and that the set $\{\ol{x}_j, \ol{y}_j \mid j \in [1,n]\}$ generates $\gr A^{E,B}_n(T)$.

The inclusion $\subseteq$ in \eqref{cent} follows from the fact that $\gr(\ZZ(A^{E,B}_n(T'))) \subseteq \ZZ(\gr A^{E,B}_n(T'))$. The
reverse inclusion is proved by using that $z_j$, $x_j^{d_j}$ and $y_j^{d_j}$ are normal elements 
of $A^{E,B}_n(T)$ which normalize the generators $x_k$, $y_k$, and then applying the first four relations in \eqref{rel} and \eqref{xy-normal}.

An isomorphism of $T'$-algebras $\ZZ(A^{E, B}_n (T')) \cong \ZZ( \gr  A^{E, B}_n (T'))$ is given by $r \mt \gr r$ where 
$r$ ranges over the $T'$-basis of $\ZZ(A^{E,B}_n(T'))$ given by the RHS of \eqref{cent}. This follows from \eqref{cent-gr}, \eqref{cent} 
and the fact that the elements $x_j^{d_j}$, $z_j$ and $y_j^{d_j}$ normalize each other.
\end{proof}
%%%%%%%%%%%
\subsection{Proof of Theorem A (ii).}
\label{3.2}
(I) The equivalence (a) $\Leftrightarrow$ (b) in Theorem A (ii) is a general fact for skew polynomial algebras \cite[Lemma 2.3]{CPWZ2}. 
%The equivalences (a) $\Leftrightarrow$ (c) and (b) $\Leftrightarrow$ (d) follow from \thref{Centers}.

(II) Next we show that (c) $\Rightarrow$ (a) $\Rightarrow$ (b) $\Rightarrow$ (e) $\Rightarrow$ (c). 
If (c) is satisfied, then  $A^{E,B}_n(T')$ is a free module over $\ZZ(A^{E,B}_n(T'))$. 
By considering the leading terms of a basis and using that $\ZZ(A^{E, B}_n (T')) \cong \ZZ( \gr  A^{E, B}_n (T'))$ (\thref{Centers}), 
we obtain that (a) holds. Thus, (c) $\Rightarrow$ (a).

By \cite[Lemma 2.3]{CPWZ2}, (a) $\Rightarrow$ (b) and if (b) is satisfied, then  
$\ZZ( \gr A^{E,B}_n(T))$ is a polynomial algebra over $T$ in powers of the generators $\ol{x}_1$, $\ol{y}_1,$ $\ldots$, 
$\ol{x}_n$, $\ol{y}_n$. So, (b) implies that
\[
\ZZ( \gr A^{E,B}_n(T)) = T[\ol{x}_1^{l'_1}, \ol{y}_1^{l_1}, \ldots, \ol{x}_n^{l'_n}, \ol{y}_n^{l_n}]
\] 
for some $l'_1, l_1, \ldots, l'_n, l_n \in \Zset_+$. It follows from \thref{Centers} that
\[
(\ol{y}_{j-1} \ol{x}_{j-1})^{c} (\ol{y}_j \ol{x}_j)^{d_j} \in \ZZ( \gr A^{E,B}_n(T))
\]
for all $c \in \Nset$ such that $d_k | (c + d_j)$ for $k <j$. Therefore $l'_j | d_j$ and $l_j | d_j$. 
However, the last commutation relation in \eqref{rel} gives that $d_j | l'_j$ and $d_j | l_j$. Thus, $l'_j = l_j =d_j$ for all 
$j \in [1,n]$ and 
\[
\ZZ( \gr A^{E,B}_n(T)) = T[\ol{x}_1^{d_1}, \ol{y}_1^{d_1}, \ldots, \ol{x}_n^{d_n}, \ol{y}_n^{d_n}].
\]
Since $\ol{x}_j^{d_j}, \ol{y}^{d_j} \in \ZZ( \gr A^{E,B}_n(T))$, from the first and fourth relations in \eqref{rel} we obtain
\[
\be_{jk}^{d_j} = 1 \; \; \mbox{for} \; \; k \in [1,n], 
\quad (\ep_k \be_{kj})^{d_j} = 1 \;  \; \mbox{for} \; \; k<j.
\]
Therefore, $\ep_k^{d_j} = 1$ for $k <j$. So $d_{jk} | d_j$ for all $k$ and $d_k | d_j$ for all $k<j$. 
This shows that (b) $\Rightarrow$ (e).

Now, assume (e). We will show that 
\begin{equation}
\label{CEB-ident}
C(E,B) = \{ (a_1, b_1, \ldots, a_n, b_n) \in \Nset^{2n} \mid \; d_j | a_j \; \; \mbox{and} \; \; d_j | b_j, \; \; \forall j \}.
\end{equation}
By the first condition in the definition of $C(E,B)$, $d_j \mid (a_j - b_j)$ and by (e), $d_{jk} \mid d_j$, $\forall k$.
Therefore the second condition in the definition of $C(E,B)$ reduces to 
\[
d_k \mid (a_k + \cdots + a_n), \; \; \forall k.
\]
Working backwards and using that $d_k \mid d_j$ for $k<j$ one gets that this condition is equivalent to: 
$d_k | a_k$ for all $k$. \thref{Centers} and \eqref{CEB-ident} imply \eqref{center}, which in turn 
implies (c). This proves that (e) $\Rightarrow$ \eqref{center} $\Rightarrow$ (c).

(III) Finally we show that (d) $\Rightarrow$ (b) and (e) $\Rightarrow$ (d). Assume that (d) is satisfied. 
It follows from \eqref{cent} that
\[
\ZZ(A^{E,B}_n(T')) \cong \ZZ(A^{E,B}_n(T)) \otimes_T T',
\]
so, $\ZZ(A^{E,B}_n(T'))$ is a polynomial algebra over $T'$. Using the isomorphism 
of $T'$-algebras $\ZZ(A^{E, B}_n (T')) \cong \ZZ( \gr  A^{E, B}_n (T'))$ (\thref{Centers}),
we obtain that $\ZZ(\gr A^{E,B}_n(T'))$ is a polynomial algebra over $T'$. The property that 
the center of a skew polynomial ring is a polynomial algebra depends on its structure 
constants and not on the base ring. This implies that $\ZZ(\gr A^{E,B}_n(T))$ is 
a polynomial algebra over $T$, which shows that (d) $\Rightarrow$ (b). 

We proved that (e) $\Rightarrow$ \eqref{center} and we have that \eqref{center} $\Rightarrow$ (d). 
Thus, (e) $\Rightarrow$ (d), which completes the proof of the theorem. 
\qed 

\bex{Nonpolyn} Consider the case $n=2$, $\be_{jk} =1$ for all $j,k$. In this case  
\begin{align*}
C(E,B) = \big( &\Zset(d_1, 0, 0, 0) + \Zset(0, d_1, 0, 0) + \\
&+ \Zset(-d_2, -d_2, d_2, 0) + \Zset(0,0,0,d_2)   \big) \cap \Nset^4.
\end{align*}
If $d_1 \nmid d_2$, then $\ZZ(\gr A^{E,B}_2(T))$ and  $\ZZ(A^{E,B}_2(T))$ are not polynomial rings. 
\eex
%%%%%%%%%%%%%
\subsection{The induced Poisson structure on $\ZZ(A^{E,B}_n(T))$}
\label{3.3}
Denote by $A^{E,B}_n(T)_q$ the algebra over $T[q^{\pm 1}]$ with 
generators $\wt{x}_j$, $\wt{y}_j$ and relations as in \eqref{rel} where $\ep_j$ are replaced by $q^{d_n m_j /d_j}$ and 
$\be_{jk}$ by $q^{d_n m_{jk}/d_{jk}}$.
(Here $q$ is an indeterminate.)
It is well known that $A^{E,B}_n(T)_q$ is a domain. 
Let $\ep := \exp( 2 \pi \sqrt{-1}/d_n) \in T^\times$. So, 
\[
\ep_j = \ep^{d_n m_j /d_j} \quad \mbox{and} \quad \be_{jk} = \ep^{d_n m_{jk}/d_{jk}}.
\]
The following lemma realizes 
$A^{E,B}_n(T)$ as the specialization of $A^{E,B}_n(T)_q$ at $q = \ep$. Its proof is straightforward and 
is left to the reader.
\ble{spe} For any integral domain $T$, there is a unique $T$-algebra homomorphism 
$\sig \colon A^{E,B}_n(T)_q \to A^{E,B}_n(T)$ such that $\sig(\wt{x}_j) = x_j$, $\sig(\wt{y}_j) = y_j$, 
$\sig(q) = \ep$. Its kernel is 
\[
\ker \sig = (q-\ep) A^{E,B}_n(T)_q.
\]
The map $\sig$ realizes $A^{E,B}_n(T)$ as the specialization of $A^{E,B}_n(T)_q$ at $q = \ep$.
\ele
Recall the definitions \eqref{XYgen}--\eqref{Zgen} of $X_j, Y_j, Z_j \in \ZZ(A^{E,B}_n(T))$. 
\bpr{Poisson} Let $T$ be an integral domain. Assume that the conditions in Theorem A (ii) 
are satisfied. Then
\begin{equation}
\label{specz}
\sig(\wt{z}_j)^{d_j} = z_j^{d_j}=Z_j \in \ZZ(A^{E,B}_n(T)).
\end{equation}
The induced Poisson bracket on $\ZZ(A^{E,B}_n(T))$ is given by
\begin{align}
\{X_j, Y_j \} &= m_j d_j d_n \ep^{-1} (X_j Y_j - (1-\ep_j)^{-d_j} Z_{j-1}), \quad &\forall j,
\nn
\\
\{Y_j, Y_k\} &=  \frac{m_{jk} d_j d_k d_n}{d_{jk}} \ep^{-1} Y_j Y_k, \quad                         &\forall j,k,
\nn
\\
\{X_j, Y_k\} & = - \frac{m_{jk} d_j d_k d_n}{d_{jk}} \ep^{-1} X_j Y_k, \quad &j<k,
\label{bracxy}
\\
\{X_j, X_k\} &=  d_j d_k d_n \left( \frac{m_j}{d_j} + \frac{m_{jk}}{d_{jk}} \right) \ep^{-1} X_j X_k, \quad &j<k,
\nn
\\
\{X_j, Y_k\} &= d_j d_k d_n \left( \frac{m_k}{d_k} + \frac{m_{kj}}{d_{kj}} \right)\ep^{-1}  X_j Y_k, \quad &j>k.
\nn
\end{align}
It satisfies
\begin{align}
\label{bracz}
\{Z_j, X_k\} &= -\de_{k \leq j} m_k d_j d_n \ep^{-1} Z_j X_k, \\
\{Z_j, Y_k\} &= \de_{k \leq j} m_k d_j d_n \ep^{-1} Z_j Y_k, \quad
\nn
\{Z_j, Z_k \} =0, \quad \forall j,k.
\end{align}
\epr
\begin{proof} Using \eqref{xyz}, we obtain 
\[
z_j^{d_j} = ((\ep_j-1) y_j x_j + z_{j-1})^{d_j} = 
(\ep_j -1)^{d_j} \ep_j^{d_j (d_j-1)/2} y_j^{d_j} x_j^{d_j}
+ \sum_{i=1}^{d_j-1} t_i y_j^i x_j^i z_{j-1}^{d_j-i} + z_{j-1}^{d_j} 
\]
for some $t_i \in T$. Since $\ZZ(A^{E,B}_n(T))$ is generated by 
$x_k^{d_k}, y_k^{d_k}$, $k \in [1,n]$ and $z_j^{d_j} \in \ZZ(A^{E,B}_n(T))$, 
$t_i=0$ for all $i$. (One can also verify this directly using $q$-binomial identities.)
Thus, 
\[
\sig(\wt{z}_j)^{d_j}= z_j^{d_j} = (\ep_j-1)^{d_j} \ep_j^{d_j(d_j-1)/2} y_j^{d_j} x_j^{d_j} + z_{j-1}^{d_j}
\]
and \eqref{specz} follows by induction on $j$ from the definition \eqref{Zgen} of $Z_j$ and the identity 
\begin{equation}
\label{pm}
\ep_j^{1 + \cdots + (d_j-1)} = (-1)^{d_j -1}.
\end{equation}

Similarly, we have 
\begin{multline*}
\wt{x}_j^{d_j} \wt{y}_j^{d_j} - q^{m_j d_j d_n} \wt{y}_j^{d_j} \wt{x}_j^{d_j} = \sum_{i=0}^{d_j-1} t_i \wt{y}_j^i \wt{x}_j^i \wt{z}_{j-1}^{d_j -i} 
\\
\mbox{for} \quad t_0 = \prod_{k=1}^{d_j}(1+ q^{d_n m_j /d_j} + \cdots + q^{(k-1)d_n m_j /d_j}), t_i \in T[q^{\pm 1}].
\end{multline*}
From here we obtain that
\[
\{X_j, Y_j \} = \sig \left( \frac{x_j^{d_j} y_j^{d_j} - y_j^{d_j} x_j^{d_j}}{q- \ep} \right)
= m_j d_j d_n \ep^{-1} (X_j Y_j - (1-\ep_j)^{-d_j} Z_{j-1})
\]
because \eqref{center} implies that $\sig(t_i \wt{y}_j^i \wt{x}_j^i \wt{z}_{j-1}^{d_j-i}  /(q-\ep)) =0$ for $0 < i < d_j$.

The rest of the Poisson brackets between any two of the functions $\{X_j, Y_j \mid j \in [1,n] \}$ follow from the defining relations of $A^{E,B}_n(T)_q$. 
The Poisson brackets in \eqref{bracz} are derived from the normalizing identities in $A^{E,B}_n(T)_q$ 
\[
\wt{z}_j \wt{x}_k = q^{-\de_{k \leq j}d_n/d_k} \wt{x}_k \wt{z}_j , \quad 
\wt{z}_j \wt{y}_k = q^{\de_{k \leq j} d_n/d_k} \wt{y}_k \wt{z}_j, \quad 
\wt{z}_j \wt{z}_k = \wt{z}_k \wt{z}_j,
\]
valid for all $j$, $k$, which are analogous to \eqref{z-normal}.
\end{proof}
%%%%%%%%%%%%%%%%%%%%%%%%%%%%%%%%%%%%
\sectionnew{Poisson prime elements of $\ZZ(A^{E,B}_n(T))$ and discriminants}
\label{discr-qWeyl}
This section contains a classification of the Poisson prime elements of the centers of the PI quantized Weyl algebras 
and our first proof of Theorem B (ii) using Poisson geometry.
\subsection{Poisson prime elements of $\ZZ(A^{E,B}_n(T))$}
\label{4.1}
Assume the conditions in Theorem A (ii) and recall the definitions \eqref{XYgen}--\eqref{Zgen} 
of the elements $X_j,Y_j,Z_j \in \ZZ(A^{E,B}_n(T))$. If $T= \Cset$, then we can consider the Poisson structure $\pi$ on 
\[
\Spec \ZZ(A^{E,B}_n(\Cset)) \cong \Aset^{2n}
\]
corresponding to the Poisson bracket from \prref{Poisson}. The Poisson brackets
\begin{equation}
\label{Pbrack}
\{X_j, X_k\} \in \Cset[X_1, \ldots, X_n], \; \; \{Z_j, Z_k \} =0, \; \; 
\{Z_j, X_k\} = -\de_{k \leq j} d_j d_k d_n Z_j X_k
\end{equation}
imply that the Poisson structure $\pi$ is symplectic on 
\begin{equation}
\label{compl}
\Aset^{2n} \backslash \big( (\cup_j \VV(X_j)) \cup (\cup_j \VV(Z_j)) \big).
\end{equation}
The set \eqref{compl} is contained in a single symplectic leaf of $\pi$, because it is connected.
By the comment after \deref{prime-norm}, every Poisson prime element of $(\ZZ(A^{E,B}_n(\Cset)), \{.,.\})$ should be an associate 
of $X_j$ or $Z_j$; however, the elements $X_j$ are not Poisson normal by \prref{Poisson}. 
Thus, the only Poisson prime elements of $\ZZ(A^{E,B}_n(\Cset))$ (up to associates) are $Z_1, \ldots, Z_n$. The idea of the proof of 
Theorem B (i) for any base field of charactecteristic 0 is to translate this argument to an algebraic setting, to avoid reference to 
Poisson geometry. This was suggested to us by Ken Goodearl.

\bpr{Poisson-prim} If $T$ is a field of characteristic 0 and the conditions in Theorem A (ii) are statisfied, 
then the following hold:

(i) The Poisson algebra 
\[
A^\circ:= \ZZ(A^{E,B}_n(T))[X_j^{-1}, Z_j^{-1}, 1 \leq j \leq n]
\]
is Poisson simple, i.e., it has no nontrivial Poisson ideals. 

(ii) Up to associates, the Poisson prime elements of $(\ZZ(A^{E,B}_n(T)), \{.,.\})$ are $Z_1, \ldots, Z_n$.
\epr
\begin{proof} The linearity of the recursion \eqref{Zgen} implies that the elements $Z_j \in \ZZ(A^{E,B}_n(\Cset))$ are irreducible and thus prime. 
Any $f \in A^\circ$, can be written in a unique way in the form
\[
f = \sum g_i h_i, 
\] 
where $g_i \in T[Z_1^{\pm 1}, \ldots, Z_n^{\pm 1}]$ and $h_i$ are distinct Laurent monomials in $X_1, \ldots, X_n$. It follows from the second and third set
of Poisson brackets  in \eqref{Pbrack} that, if $f$ belongs to a Poisson ideal $I$ of $A^\circ$, then all terms $g_i h_i$ belong to $I$ and, thus, 
$g_i \in I$, $\forall i$.
Therefore, if $I$ is a nonzero Poisson ideal of $A^\circ$, then $I \cap \Cset[Z_1^{\pm1}, \ldots, Z_n^{\pm 1}] \neq \{ 0 \}$. 
Exchanging the roles of $X_j$ and $Z_j$ and using the first and third Poisson brackets in \eqref{Pbrack} implies that $I = A^\circ$ which 
proves part (i). 

(ii): From the second and third set of Poisson brackets in \eqref{Pbrack} we obtain that $Z_1, \ldots, Z_n$ are Poisson prime elements of $(\ZZ(A^{E,B}_n(T)), \{.,.\})$.
\prref{Poisson} implies that the elements $X_j$ are not Poisson normal, and thus they are not Poisson prime. 
If $f \in \ZZ(A^{E,B}_n(T))$ is a Poisson prime element that is not an associate of any $Z_j$, then $f$ would also be a Poisson 
prime element of $A^\circ$ because $f \nmid X_1 \ldots X_n Z_1 \ldots Z_n$. 
This means that $(f)$ would be a nontrivial Poisson ideal of $A^\circ$ which contradicts part (i).
\end{proof}
%%%%%%%%%%%
\subsection{First proof of Theorem B (ii)}
\label{4.2}
Denote by $\tr \colon A^{E,B}_n(T) \to \ZZ(A^{E,B}_n(T))$ the internal trace associated to the embeddings $A^{E,B}_n(T) \hra M_{N^2} (\ZZ(A^{E,B}_n(T)))$ 
from $\ZZ(A^{E,B}_n(T))$-bases of $A^{E,B}_n(T)$, recall that $N= d_1 \ldots d_n$. The set
\[
\B := \{ x_1^{l'_1} y_1^{l_1} \ldots x_n^{l'_n} y_n^{l_n} \mid l_j, l'_j \in [0,d_j-1] \} 
\]
is a $\ZZ(A^{E,B}_n(T))$-basis of $A^{E,B}_n(T)$. We will prove that 
\begin{align}
\ep^m d_{N^2}( \B : \tr) &= \eta Z_1^{N^2 (d_1-1)/d_1} \ldots Z_n^{N^2 (d_n-1)/d_n}
\label{dis}
\\
&=  \eta z_1^{N^2 (d_1-1)} \ldots z_n^{N^2 (d_n-1)} 
\nn
\end{align}
for some $m \in \Zset$ (depending only on $E$ and $B$), where $\eta \in T$ is the scalar defined in Theorem B (ii) and 
$\ep = \exp( 2 \pi \sqrt{-1}/d_n) \in T^\times$ as in \S \ref{3.3}. We note that 
\begin{equation}
\label{eta2}
\eta = \Big( N^2  (1 - \ep_1)^{-d_1+1} \dots (1- \ep_n)^{-d_n+1} \Big)^{N^2},
\end{equation}
because 
\begin{equation}
\label{ep1}
(1- \ep_j)^{d_j-1} \prod_{i=1}^{d_j-2} (1 + \ep_j + \cdots + \ep_j^i)  = \prod_{i=1}^{d_j-1} (1 - \ep_j^i) = d_j.
\end{equation}
In particular, $d_j (1-\ep_j)^{-d_j+1} \in T^\times$.
Since $A^{E,B}_n(T)$ is defined over $\Zset[\ep]$ and $\B \subset A^{E,B}_n(\Zset[\ep])$, 
it is sufficient to prove \eqref{dis} for $T = \Zset[\ep]$. From now on we will assume that $T = \Zset[\ep]$.

Recall the filtration of $A^{E,B}_n(T)$ from \S \ref{1.1}. 
By \cite[Proposition 4.10]{CPWZ2}, 
\[
\gr d_{N^2} (\B : \tr) = d_{N^2} (\gr \B : \tr).
\]
Note that $\gr \B$ is a $\ZZ( \gr A^{E,B}_n(T))$-basis of $\gr A^{E,B}_n(T)$ by Theorem A (i); the RHS uses the internal trace of $\gr A^{E,B}_n(T)$
with respect to this basis.
Since $\gr A^{E,B}_n(T)$ is a skew polynomial algebra, we can apply \cite[Proposition 2.8]{CPWZ2} to deduce that
\[
\ep^{m'} d_{N^2} (\gr \B : \tr) = N^{2N^2} (\ol{x}_1 \ol{y}_1)^{N^2 (d_1-1)} \ldots (\ol{x}_n \ol{y}_n)^{N^2 (d_n-1)}
\]
for some $m' \in \Zset$. We have $\gr z_j = (\ep_j -1) \ol{x}_j \ol{y}_j$. It follows from \eqref{eta2} that
\begin{equation}
\label{gr-rel}
\ep^{m'} 
\gr d_{N^2} (\B : \tr) = \pm \eta \gr (z_1^{N ^2(d_1-1)} \ldots z_n^{N^2 (d_n-1)}),
\end{equation} 
where $\pm = (-1)^{N^2(d_1 + \cdots + d_n -n)}$, which is a power of $\ep$ by \eqref{pm}.
By passing from $T= \Zset[\ep]$ to its field of fractions, we see that to prove Theorem B (ii), it is sufficient to prove it in the case when 
$T= \Qset(\ep)$. From now on we will assume that $T= \Qset(\ep)$.
Combining \thref{prod} and \prref{Poisson-prim} (ii), we obtain
\begin{equation}
\label{inter}
d(A^{E,B}_n(\Cset)/ \ZZ(A^{E,B}_n(\Cset))) =_{\Qset(\ep)^\times} Z_1^{s_1} \ldots Z_n^{s_n}
\end{equation}
for some $s_1, \ldots, s_n \in \Nset$. \prref{Poisson} and eq. \eqref{gr-rel} imply that $s_j = N^2 (d_j -1)/d_j$. The first equality 
in \eqref{dis} for $T= \Qset(\ep)$ follows from \eqref{gr-rel} and \eqref{inter}. By the above argument this fact holds for all 
integral domains $T$. The second equality in \eqref{dis} 
follows from the fact that $Z_j = z_j^{d_j}$, recalling \prref{Poisson}. \qed

\bre{twist} One can use 2-cocycle twists \cite{AST} to give a conceptual proof of the fact that 
the discriminant formula in Theorem B (ii) does not depend on the matrix $B$. Denote  by 
$1_n$ the $n \times n$ matrix, all entries of which equal $1$. The algebras $A^{E,1_n}_n(T)$ and 
$A^{E, B}_n(T)$ are $\Zset^n$-graded by 
\[
\deg y_j =  - \deg x_j = e_j,
\] 
where $e_1, \ldots, e_n$ is the standard basis of $\Zset^n$. The algebra $A^{E,B}_n(T)$ is obtained 
from $A^{E,1_n}_n(T)$ by the 2-cocycle twist \cite{AST} using the cocycle
\[
\ga \colon \Zset \times \Zset \to T^{\times}, \quad \ga(e_j, e_k) = \ga(e_k, e_j)^{-1} := \sqrt{\be_{jk}}, \; \; j <k. 
\]
If $A^{E,B}_n(T)$ satisfies the conditions of Theorem A (ii), then the same is true for $A^{E,1_n}_n(T)$. 
The centers of both algebras are generated by $x_j^{d_j}$, $y_j^{d_j}$ and the twist on $\ZZ^{E,1_n}_n(T)$ 
is trivial, meaning the product is not changed under the twist. One easily checks that for degree reasons, the 
trace of a homogeneous element of $A^{E,1_n}_n(T)$ is the same under the two products 
(the one in $A^{E,1_n}_n(T)$ and the twisted one). This implies that the two discriminants  
$d(A^{E,B}_n(T)/ \ZZ(A^{E,B}_n(T)))$ and $d(A^{E,1_n}(T)/ \ZZ(A^{E,1_n}_n(T)))$ 
are equal.
A similar argument can be given for the independence of the discriminant formula in \thref{general-discr} on the 
entries of the matrix $B$. 
\ere
%%%%%%%%%%%%%%%%%%%%%%%%%%%%%%%%%%%%
\sectionnew{Discriminants of quantized Weyl algebras using quantum cluster algebras}
\label{q-cl}
%%%%%%%%
In this section we give a second proof of Theorem B (ii) using techniques from quantum cluster algebras.
%%%%%%%%%
\subsection{Quantum cluster algebras and quantized Weyl algebras}
\label{4.1a}
Cluster algebras were introduced by Fomin and Zelevinsky in \cite{FZ}. Their quantum counterparts were defined
by Berenstein and Zelevinsky 
in \cite{BeZ}. A quantum cluster algebra has (generally infinitely many) localizations that are isomorphic 
to quantum tori. The generators of the different quantum tori are related by mutations. The set of all 
such generators of quantum tori generate the whole quantum cluster algebra.

Our second proof of Theorem B (ii) relies on relating the discriminant of a quantum cluster algebra to the discriminants of the 
corresponding quantum tori and then evaluating the former, the key point being that the later are straightforward to compute. 
The outcome of this is a formula for the discriminant in terms of a product of frozen cluster variables. 
   
When the parameters $\ep_j$ are non-roots of unity, the quantized Weyl algebras $A^{E,B}_n(T)$ are symmetric CGL extensions 
and the result in \cite[Theorem 8.2]{GY2} constructs quantum cluster algebra structures on $A^{E,B}_n(T)$ when $T$ is a field. We will work with the 
PI case. We will only use two quantum clusters whose intersection is precisely the set of frozen variables.

\bpr{qcl-Weyl} Let $T$ be an integral domain, $n \in \Zset_+$, and $\ep_j, \be_{jk} \in T^\times$, $\ep_j \neq 1$ be such that $\ep_j -1 \in T^\times$
for all $j$.
\begin{enumerate}
\item[(i)]
The localization of $A^{E,B}_n(T)[y_j^{-1}, 1 \leq j \leq n]$ is isomorphic to 
the mixed skew-polynomial/quantum torus algebra over $T$ with generators 
\[
y_j^{\pm1}, z_j, \quad j \in [1,n]
\]
and relations
\begin{equation}
\label{y-zrel}
y_j y_k = \be_{jk} y_k y_j, \; \; z_j z_k = z_k z_j, \; \; z_j y_k = \ep_k^{\de_{k \leq j}} y_k z_j, \; \; j, k \in [1,n].
\end{equation}
\item[(ii)] The localization of $A^{E,B}_n(T)[x_j^{-1}, 1 \leq j \leq n]$ is isomorphic to 
the mixed skew-polynomial/quantum torus algebra over $T$ with generators $x_j^{\pm1}$, $z_j$, $j \in [1,n]$
and relations
\[
x_j x_k = \ep_j \be_{jk} x_k x_j, \; \; j< k; \; \; z_j z_k = z_k z_j, \; \; z_j x_k = \ep_k^{-\de_{k \leq j}} x_k z_j, \; \; 
j, k \in [1,n].
\]
\end{enumerate}
\epr
\begin{proof} (i) It follows from the definition of $A^{E,B}_n(T)$ and \eqref{z-normal} that the elements 
$y_j, z_j \in A^{E,B}_n(T)$ satisfy the stated relations. The isomorphism follows from the fact that the generators
$x_j$ of $A^{E,B}_n(T)$ can be expressed in terms of the elements $y_j^{\pm1}, z_j$ by
\[
x_j = (\ep_j-1)^{-1} y_j^{-1} ( z_j - z_{j-1}). 
\]
Part (ii) is analogous.
\end{proof}
The quantum clusters in parts (i) and (ii) of \prref{qcl-Weyl} correspond to the ones constructed in 
\cite[Theorem 1.2]{GY1} from the CGL extension presentations of $A^{E,B}_n(T)$ associated to adjoining 
its generators in the orders 
\[
y_n, \ldots, y_1, x_1, \ldots, x_n \quad \mbox{and} \quad x_n, \ldots, x_1, y_1, \ldots, y_n, 
\]
respectively. Technically, \prref{qcl-Weyl} allows the scalars $\ep_1, \ldots, \ep_n$ to be roots of unity 
different from 1, while the general result in \cite{GY1} requires those to be non-roots of unity.
%%%%%%%%%
\subsection{A second proof of Theorem B (ii)}
\label{4.2a}
Similarly to \S \ref{4.2}, one shows that it is sufficient to prove the 
theorem in the case when $\ep_j -1 \in T^\times$ for all $j$. Indeed, start from an arbitrary 
integral domain $T$, denote
\[
T' := T[(\ep_j -1)^{-1}, 1 \leq j \leq n],
\]
and assume that the theorem is valid when the base ring is $T'$. Because of \eqref{eta2}, this means that
\[
d(A^{E,B}_n(T')/\ZZ(A^{E,B}_n(T'))) =_{(T')^\times} N^{2 N^2} z_1^{N^2 (d_1-1)} \ldots z_n^{N^2 (d_n-1)}.
\]
Therefore, 
\[
d(A^{E,B}_n(T)/\ZZ(A^{E,B}_n(T))) =_{T^\times} \nu N^{2 N^2} z_1^{N^2 (d_1-1)} \ldots z_n^{N^2 (d_n-1)}.
\]
for some $\nu \in T'$. By \cite[Propositions 2.8 and 4.10]{CPWZ1} (same argument as the one for \eqref{gr-rel}),
\[
\gr d(A^{E,B}_n(T)/\ZZ(A^{E,B}_n(T))) =_{T^\times} \eta (\gr z_1)^{N^2 (d_1-1)} \ldots (\gr z_n)^{N^2 (d_n-1)}
\]
with respect to the filtration from \S \ref{1.1}. This implies that $\nu N^{2 N^2} =_{T^\times} \eta$ and 
\[
d(A^{E,B}_n(T)/\ZZ(A^{E,B}_n(T))) =_{T^\times} \eta z_1^{N^2 (d_1-1)} \ldots z_n^{N^2 (d_n-1)}.
\]

{\em{From now on we will assume that $\ep_j - 1 \in T^\times$ for all $j$.}}
%By \eqref{ep1}, this is equivalent to $N = d_1 \ldots d_n \in T^\times$, so $\eta \in T^\times$. 
Denote by $A(y,z,T)$ the skew polynomial algebra over $T$ with generators $y_j, z_j$, $j \in [1,n]$ and relations \eqref{y-zrel}. 
\prref{qcl-Weyl} (i) implies that  
\[
A^{E,B}_n(T)[y_j^{-d_j}, 1 \leq j \leq n] \cong A(y,z,T) [y_j^{-d_j}, 1 \leq j \leq n].
\]
Denote by $A$ this algebra. Since $\ZZ(A^{E,B}_n(T)) = T[ x_j^{d_j}, y_j^{d_j}, 1 \leq j \leq n]$, it follows from this, \eqref{Zgen} and \eqref{specz} that 
\[
\ZZ(A) = T[y_j^{\pm d_j}, z_j^{d_j}, 1 \leq j \leq n] 
\quad 
\mbox{and}
\quad
\ZZ(A(y,z,T)) = T[y_j^{d_j}, z_j^{d_j}, 1 \leq j \leq n].
\]
The algebras $A$ and $A(y,z,T)$ are free over their centers. Because $A$ is a 
central localization of $A^{E,B}_n(T)$ and $A(y,z,T)$, the internal trace
$\tr \colon A \to \ZZ(A)$ is an extension of the internal traces
$\tr \colon A^{E,B}_n(T) \to \ZZ(A^{E,B}_n(T))$ and 
$\tr \colon A(y,z, T) \to \ZZ(A(y,z,T))$. Moreover,
\begin{equation}
\label{3discr}
d(A^{E,B}_n(T)/\ZZ(A^{E,B}_n(T))) =_{\ZZ(A)^\times} 
d(A/\ZZ(A)) =_{\ZZ(A)^\times} d(A(y,z,T)/ \ZZ(A(y,z,T))).
\end{equation}
By \cite[Proposition 2.8]{CPWZ2}
\[
d(A(y,z,T)/ \ZZ(A(y,z,T)))
 =_{T^\times} N^{2 N^2} z_1^{N^2 (d_1-1)} \ldots z_n^{N^2 (d_n-1)},
\]
recall that $N \in T^\times$. It follows from \eqref{3discr} that
\[
d(A^{E,B}_n(T)/\ZZ(A^{E,B}_n(T))) =_{T^\times} N^{2 N^2} y_1^{h_1 d_1} \ldots y_n^{h_n d_n} z_1^{N^2 (d_1-1)} \ldots z_n^{N^2 (d_n-1)}
\]
for some $h_j \in \Zset$. 
Analogously, \prref{qcl-Weyl} (ii) implies that 
\[
d(A^{E,B}_n(T)/\ZZ(A^{E,B}_n(T))) =_{T^\times} N^{2 N^2} x_1^{g_1 d_1} \ldots x_n^{g_n d_n} z_1^{N^2 (d_1-1)} \ldots z_n^{N^2 (d_n-1)}
\]
for some $g_j \in \Zset$. Since $A^{E,B}_n(T)$ is a domain, 
\[
y_1^{h_1 d_1} \ldots y_n^{h_n d_n} =_{T^\times} x_1^{g_1 d_1} \ldots x_n^{g_n d_n}.
\]
This is only possible if $g_j = h_j=0$  for all $j$ because
\[
\{ y_1^{l_1} \ldots y_n^{l_n} x_1^{l'_1} \ldots x_n^{l'_n} \mid l_j, l'_j \in \Nset \}
\]
is a $T$-basis of $A^{E,B}_n(T)$ (a consequence of the fact that the latter is an iterated skew polynomial extension).
Thus, 
\[
d(A^{E,B}_n(T)/\ZZ(A^{E,B}_n(T))) =_{T^\times} N^{2 N^2} z_1^{N^2 (d_1-1)} \ldots z_n^{N^2 (d_n-1)},
\]
which proves Theorem B (ii).
\qed
\bre{q-cl-other} The proof of Chan, Young and Zhang \cite{CYZ} of the case $n=1$ of Theorem B (ii) also implicitly used
quantum cluster algebras (though the two proofs appear to be quite different). They worked with the elements $y_1$ and $y_1^{-1} z_1$ 
which generate a quantum torus. One of the clusters in the above proof consists of the cluster variables $y_1,z_1$. The 
quantum tori from the two proofs are the same.
\ere
%%%%%%%%%%%%%%%%%%%%%%%%%%%%%%%%%%%%
\sectionnew {Automorphisms and Isomorphisms between PI quantized Weyl algebras}
In this section we prove Theorems B (iii) and C.
%%%%%%%%
\subsection{Local dominance of discriminants of PI quantized Weyl algebras}
\label{5.1}
The following proposition proves Theorem B (iii).
\bpr{preprt}
Let $A^{E_1,B_1}_{n_1}(T), \ldots, A^{E_l,B_l}_{n_l}(T)$ be a set of quantized Weyl algebras over an integral domain $T$ of characteristic 0, 
satisfying the conditions in Theorem A (ii). Let $A$ be their tensor product over $T$. The discriminant $d(A, \ZZ(A))$ is 
\begin{enumerate}
\item[(i)] locally dominating and
\item[(ii)] effective.
\end{enumerate}
\epr
Denote by $z(A) \subset A$ the union of the collection of elements $z_j$ for all algebras  $A^{E_1,B_1}_{n_1}(T), \ldots, A^{E_l,B_l}_{n_l}(T)$.
Similarly, denote by $Z(A) \subset \ZZ(A)$ the union of the collection of elements $Z_j$ for these algebras. Let $x(A) \subset A$ and $y(A) \subset A$ 
denote the collections of all $x$ and $y$-generators of $A$. Since
\begin{equation}
\label{tensprod}
d(A/\ZZ(A)) = d(A^{E_1,B_1}_{n_1}(T)/ \ZZ(A^{E_1,B_1}_{n_1}(T))) \ldots d(A^{E_l,B_l}_{n_l}(T), \ZZ(A^{E_l,B_l}_{n_l}(T))),
\end{equation}
\prref{preprt} (i) directly follows from the following lemma and Theorem B (ii).
\ble{preprt1} Let $\phi \in Aut_T(A)$. 

(i) For each $z \in z(A)$, $\deg \phi(z) \geq 2$. 

(ii) If $\deg \phi(x) >1$ for at least one $x \in x(A)$ or $\deg \phi(y) >1$ for at least one $y \in y(A)$, then $\deg \phi(z) >2$ for 
some $z \in z(A)$.
\ele

Denote by $\FF_j A$ the $\Nset$-filtration of $A$ defined by 
\begin{equation}
\label{filt2}
\deg x = \deg y = 1 \quad \mbox{for all} \; \; x \in x(A), y \in y(A).
\end{equation} 
Note that this is a different filtration from the one defined in \S \ref{1.1}.

\begin{proof}[Proof of \leref{preprt1}] Part (i): It is straightforward to verify that the only normal elements of 
$A^{E_i,B_i}_{n_i}(T)$ in $\FF_1 A^{E_i,B_i}_{n_i}(T)$ are $T. 1$ and they are central. At the same time 
$\phi(z)$ is a normal element of $A$ for every $z \in z(A)$ and $\phi(z)$ is not central. Therefore 
$\phi(z) \notin \FF_1 A$.

Part (ii): Assume that $\phi \in \Aut_T(A)$ is such that $\deg \phi(x) >1$ for at least one $x \in x(A)$ or $\deg \phi(y) >1$ for at least one $y \in y(A)$.
Then there exist $i \in [1,l]$ and $j \in [1,n_i]$ such that the $x$- and $y$-generators of $A^{E_i,B_i}_{n_i}(T)$, to be denoted 
by $x_1, y_1, \ldots, x_{n_i}, y_{n_i}$, satisfy
\[
\deg \phi(x_k) = \deg \phi(y_k) =1 \; \; \mbox{for} \; \; k<j \quad \mbox{and} \quad
\deg \phi(x_j) >1 \; \; \mbox{or} \; \; \deg \phi(y_j) >1.
\]
For $z_j := 1 + (\ep_1 -1) y_1 x_1 + \cdots + (\ep_j -1) y_j x_j$ we have, $\deg \phi(x_j) \phi(y_j) >2$ and 
\[
\deg (1 + (\ep_1 -1) \phi(y_1) \phi(x_1) + \cdots + (\ep_{j-1} -1) \phi(y_{j-1}) \phi(x_{j-1})) \leq 2.
\]
Thus $\deg \phi(z_j)>2$. 
\end{proof}

\bre{dom} The discriminants, considered in \cite{CPWZ1,CPWZ2,CYZ}, possess the property that they 
have unique leading terms in certain generating sets, i.e., they are linear combinations of monomials and the powers of all monomials
are componentwise less than those of a leading one. This property implies the (global) dominance of those discriminants by \cite[Lemma 2.2 (1)]{CPWZ1}. 
The discriminants in Theorem B (ii) do not possess this stronger property (except when $n=1$) and the proof of their local dominance 
is more involved.
\ere
%%%%%%%%
\subsection{Effectiveness of discriminants of PI quantized Weyl algebras}
\label{5.2}
Consider the trivial filtration on $A$, $F_0 A:= A$ (which is different from the filtrations in \S \ref{1.1}  and \ref{5.1}).
Let $R$ be any ``testing'' filtered PI $T$-algebra. Choose elements $\theta(x), \theta (y) \in R$ for all $x \in x(A)$, $y \in y(A)$
such that for at least one $x$ or $y$, $\theta(x) \notin F_0 R$ or $\theta(y) \notin F_0 R$. 
As in the previous subsection, there exist $i \in [1,l]$ and $j \in [1,n_i]$ 
such that the $x$- and $y$-generators $x_1, y_1, \ldots, x_{n_i}, y_{n_i}$ of $A^{E_i,B_i}_{n_i}(T)$ satisfy
\[
\deg \theta(x_k) = \deg \theta(y_k) =0 \; \; \mbox{for} \; \; k<j \quad \mbox{and} \quad
\deg \theta(x_j) >0 \; \; \mbox{or} \; \; \deg \theta(y_j) >0.
\]
Therefore,
\[
\deg (1 + (\ep_1 -1) \theta(y_1) \theta(x_1) + \cdots + (\ep_{j-1} -1) \theta(y_{j-1}) \theta(x_{j-1})) \in F_0 R, \; \; 
\theta(x_j) \theta(y_j) \notin F_0 R. 
\]
So,
\[
1 + (\ep_1 -1) \theta({y_1}) \theta({x_1}) + \cdots + (\ep_j -1) \theta({y_j}) \theta({x_j}) \notin F_0 R.
\] 
The effectiveness of the discriminant $d(A/\ZZ(A))$ follows from \eqref{tensprod} and Theorem B (ii).
%%%%%%%%
\subsection{Classification of automorphisms of tensor products of quantum Weyl algebras}
\label{5.3}
Let $A$ be a tensor product of quantized Weyl algebras as in the previous two subsections and Theorem C. Denote
\[
E(A) := \{ \ep \in T^\times \mid \ep^{\pm1} \in E_1 \cup \dots \cup E_l \}.
\]
For $r \in A$ and $\ep \in E(A) \cup \{1\}$ denote 
\begin{align*}
\LL_\ep(r) &:= \{ r' \in \FF_1 A \mid r r' = \ep r' r \},
\\
\LL(r) &:= \bigoplus_{ \ep \in E(A) \cup \{ 1 \}} \LL_\ep(r), \; \; \LL^*(r) := \bigoplus_{ \ep \in E(A)} \LL_\ep(r).
\end{align*}
\noindent
{\em{Proof of Theorem C.}} Part (ii) is straightforward to verify. 

Part (i): Let $\phi \in \Aut_T(A)$. Let $\KK$ be the field of fractions of $T$. We extend $\phi$ to a $\KK$-linear 
automorphism of $A_\KK = A \otimes_T \KK$, to be denoted by the same letter. It induces a $\KK$-linear 
automorphism of $\ZZ(A_\KK)$ which is a polynomial algebra. In addition, $\phi(d(A_\KK/\ZZ(A_\KK)) =_{\KK^\times} d(A_\KK/ \ZZ(A_\KK))$. 
The set of prime divisors of $d(A_\KK/\ZZ(A_\KK)) \in \ZZ(A_\KK)$ is $Z(A)$. Therefore, for every $Z \in Z(A)$ there exists 
$\al_0 \in \KK^\times$ such that $\al_0^{-1} \phi(Z) \in Z(A)$. 

If $z \in z(A)$, then $\phi(z)^k \in Z(A)$ for some $k \in \Zset_+$. The above implies that for every $z \in z(A)$, there exists 
$z' \in z(A)$, $k, k' \in \Zset_+$ and $\al_0 \in \KK^\times$, such that
\[
\phi(z)^k = \al_0 (z')^{k'}.
\]
Recall the filtration \eqref{filt2} of $A$.
Since $d(A_\KK/\ZZ(A_\KK))$ is locally dominating, \thref{appl-dom} implies that $\phi(\FF_1 A_\KK) = \FF_1 A_\KK$ and, thus,
$\phi(z) \in \FF_2 A_\KK$. By \leref{preprt1}, $\phi(z) \notin \FF_1 A_\KK$. 
We also have that $z' \in \FF_2 A_\KK$ and $z' \notin \FF_1 A_\KK$. Therefore, $k'= k$. Since both
$z', \phi(z) \in \FF_2 A$ are normal elements, it follows from \eqref{z-normal} that $z'$ and $\phi(z)$ commute. 
Now $\phi(z)^k = \al_0(z')^k$ implies that $z' = \al \phi(z)$ for some $\al \in \ol{\KK}$ which would have to be in $\KK$ 
since $z' \in A_\KK$. This shows that 
\begin{equation}
\label{z-invar}
\mbox{for every $z \in z(A)$ there exists $\al \in \KK^\times$ such that $\al^{-1} \phi(z) \in z(A)$.}
\end{equation} 

For each $i \in [1,l]$ there exists $z \in z(A)$ such that $\LL^*(z) = A^{E_i,B_i}_{n_i}$ and for all $z_0 \in z(A) \cap A^{E_i,B_i}_{n_i}$, 
$\LL^*(z_0) \subseteq A^{E_i,B_i}_{n_i}$. Now part (i) follows from $\phi(\FF_1 A) = \FF_1 A$ (by \thref{appl-dom} and the fact that $d(A_\KK/\ZZ(A_\KK))$ is locally dominating).

Denote by $z_1, \ldots, z_n$ and $z'_1, \ldots, z'_n$ the sequences of normal elements of $A^{E_i,B_i}_{n_i}$ 
and $A^{E_{\sig(i)}, B_{\sig(i)}}_{n_{\sig(i)}}$ that belong to $z(A)$. Set $z_0=1$. By \eqref{z-normal}, $\LL^*(z_{j-1}) \subsetneq \LL^*(z_j)$, 
$\forall j \in [1,n]$. Combining this with \eqref{z-invar} and the fact that $\phi(\FF_1 A) = \FF_1 A$, we obtain that there exist 
$\al_1, \ldots, \al_n \in \KK^\times$ such that
\begin{equation}
\label{z-invar-2}
\phi(z_j) = \al_j z'_j \quad \mbox{for} \; \; j \in [1,n].
\end{equation}

It follows from \eqref{z-normal} that for all $\ep \neq 1$
\[
\LL_\ep(z_j) \cap \LL_1(z_{j-1}) \neq 0 \; \; 
\mbox{if and only} \; \; \ep = \ep_j^{\pm1}
\]
and
\[
\LL_{\ep_j^{-1}}(z_j) \cap \LL_1(z_{j-1}) = T x_j, \quad \LL_{\ep_j}(z_j) \cap \LL_1(z_{j-1}) = T y_j.
\]
Eq. \eqref{z-invar-2} and $\phi(\FF_1 A) = \FF_1 A$ imply that either
\begin{align}
\label{case1}
&\phi(x_j)= \mu_j x_j, \quad \phi(y_j)= \nu_j y_j \quad \mbox{or} \\
&\phi(x_j)= \mu_j y_j, \quad \phi(y_j)= \nu_j x_j
\label{case2}
\end{align}
for some $\mu_j, \nu_j \in T$. The same statement for $\phi^{-1}$ gives that $\mu_j, \nu_j \in T^\times$. Set $\tau_j = 1$ in the first case 
and $\tau_j=-1$ in the second.

From the identity $[x_j, y_j] = z_j$ we obtain that $\al_j = \tau_j \mu_j \nu_j$; that is
\[
\phi(z_j) = \tau_j \mu_j \nu_j z'_j.
\]
Applying $\phi$ to both sides of the identity $x_j y_j -\ep_j y_j x_j = z_{j-1}$ and using this and \eqref{case1}-\eqref{case2} 
gives the two equalities in \eqref{ident1}. The equality \eqref{ident2} follows by applying $\phi$ to the homogeneous defining relations 
of $A^{E,B}_n(T)$ and using \eqref{case1}-\eqref{case2}.

Part (iii) follows from Theorems B (iii) and \ref{tappl-dom} (ii) and the fact that the quantized Weyl algebras have finite GK dimension.
\qed
%%%%%%%%%%%%%
\subsection{Special cases of Theorem C (i)-(ii)}
\label{5.4}
Theorem C (i)--(ii) has the following special cases classifying automorphisms and isomorphisms 
between PI quantized Weyl algebras.
\bco{1} Let $A^{E,B}_n(T)$ and $A^{E',B'}_{n'}(T)$ be two quantized Weyl algebras over an integral domain $T$ satisfying the conditions 
in Theorem A (ii), where $E = (\ep_1, \ldots, \ep_n)$, $E' = (\ep'_1, \ldots, \ep'_{n'})$, 
$B=(\be_{jk})$ and $B' = (\be'_{jk})$. The algebras $A^{E,B}_n(T)$ and $A^{E',B'}_{n'}(T)$ are isomorphic if and only if $n'=n$ and 
there exists a sequence $(\tau_1, \ldots, \tau_n) \in \{ \pm 1\}^n$ such that
\[
\ep'_j = \ep_j^{\tau_j}, \; \; \forall j \quad \mbox{and} \quad
\be'_{jk} = 
\begin{cases}
\be_{jk}^{\tau_j}, &\mbox{if} \; \; \tau_k =1
\\
(\ep_j \be_{jk})^{- \tau_j}, & \mbox{if} \; \; \tau_k =-1,
\end{cases}
\quad\quad \forall j<k.
\]
\eco
The analog of the theorem in the non-PI case was obtained in \cite{GH}. The case of the theorem when $n=1$ 
was obtained in \cite{Ga,CYZ}. The homogenized PI quantized Weyl algebras were treated in \cite{Ga2}
using the result of \cite{BZ0} on the isomorphism problem for $\Nset$ graded algebras. We note that the latter result is 
not applicable to quantized Weyl algebras because they are not $\Nset$-graded (in a nontrivial way).

\bco{2}  Let $A^{E,B}_n(T)$ be a quantized Weyl algebra over an integral domain $T$ satisfying the conditions 
in Theorem A (ii).
\begin{enumerate} 
\item For all scalars $\mu_1, \nu_1, \ldots, \mu_n, \nu_n \in T^\times$ such that $\mu_j \nu_j =1$, $\forall j$,
\[
\phi(x_j)= \mu_j x_j, \quad \phi(y_j)= \nu_j y_j
\]
defines a $T$-linear automorphism of $A^{E,B}_n(T)$. 
\item Assume that for some $k \in [1,n]$, $\ep_k=-1$, $\be_{jk}^{2} = \ep_j$ for $j < k$, $\be_{jk}^2=1$ for $j>k$. 
For all scalars $\mu_1, \nu_1, \ldots, \mu_n, \nu_n \in T^\times$ such that $\mu_j \nu_j =1$ for $j\leq k$ and $\mu_j \nu_j = -1$ for $j > k$, 
\begin{align*}
&\phi(x_j)= \mu_j x_j, \quad \phi(y_j)= \nu_j y_j, \; \; \mbox{for all} \; \; j \neq k, \\
&\phi(x_k)= \mu_k y_k, \quad \phi(y_k)= \nu_k x_k
\end{align*}
defines a $T$-linear automorphism of $A^{E,B}_n(T)$.
\end{enumerate}
All elements of $\Aut_T(A^{E,B}_n(T))$ have one of the above two forms.

In particular, $\Aut_T(A^{E,B}_n(T)) \cong (T^\times)^n \ltimes \Zset_2$ if the pair $(E,B)$ satisfies the condition in {\em{(2)}} for some $k \in [1,n]$ and 
$\Aut_T(A^{E,B}_n(T)) \cong (T^\times)^n$ otherwise.  
\eco
Note that, if the condition in (2) is satisfied, then the condition in Theorem A (ii) that $d_j | d_k$ for $j <k$ implies that $\ep_j =-1$ for all $j<k$.

The analog of the theorem in the non-PI case was obtained in \cite{R}. The case of the theorem when $n=1$ 
was obtained in \cite{CYZ}.
%%%%%%%%
\sectionnew{General discriminant formula for quantized Weyl algebras}
\label{general}
In this section we prove a general formula for the discriminants of PI quantized Weyl algebras over polynomial 
central subalgebras generated by powers of the standard generators of the Weyl algebra. 
The results are obtained by extending the approach from Sect. \ref{q-cl}, based on quantum cluster algebra techniques, combined 
with field theory.
\subsection{Formulation of main result}
\label{A.1}
To have a setting that is suitable for induction, we work with slightly more general algebras than the quantized 
Weyl algebras $A^{E,B}_n(T)$. We choose an indeterminate $c$ and 
define $A^{E,B,c}_n(T)$ to be the $T[c]$-algebra with generators $x_1, y_1, \ldots, x_n, y_n$ and relations as in \eqref{rel} 
except for the last one, which is replaced by 
\[
x_j y_j - \ep_j y_j x_j = c + \sum_{i=1}^{j-1} (\ep_i -1) y_i x_i, \quad \forall i.
\]
Here $T$ is an integral domain, as before.
The $T$-algebra $A^{E,B}_n(T)$ is obtained from it by specialization:
\[
A^{E,B}_n(T) \cong A^{E,B,c}_n(T)/(c-1) A^{E,B,c}_n(T).
\]
Assuming \eqref{scalars}, for $j <k$, represent 
\[
\frac{m_j}{d_j} + \frac{m_{jk}}{d_{jk}} = \frac{m'_{jk}}{d'_{jk}} 
\]
with $m'_{jk} \in \Nset$, $d'_{jk} \in \Zset_+$ such that $\gcd(m'_{jk}, d'_{jk}) =1$. Set $d'_{kj} := d'_{jk}$. 
Theorems A and \ref{tCenters} remain valid in this slightly more general situation.
The following lemma follows directly from Theorem A (i).
\ble{cent} For $l \in \Zset_+$ the following hold:
\begin{enumerate}
\item[(i)] $x_j^l \in \ZZ(A^{E,B,c}_n(T))$ if and only if $x_j^l \in \ZZ(A^{E,B}_n(T))$ if and only if
\[
\lcm(d_j, d'_{jk}, 1 \leq k \leq n, k \neq j) | l.
\]
\item[(ii)] $y_j^l \in \ZZ(A^{E,B,c}_n(T))$ if and only if $y_j^l \in \ZZ(A^{E,B}_n(T))$ if and only if
\[
\lcm(d_j, d_{jk}, 1 \leq k \leq n, k \neq j) | l.
\]
\end{enumerate}
\ele
For every polynomial central subalgebra of $A^{E,B,c}_n(T)$ of the form 
\begin{equation}
\label{C}
\CC := T[c,x_1^{L_1}, y_1^{L_1}, \ldots, x_n^{L_n}, y_n^{L_n}],
\end{equation}
$A^{E,B,c}_n(T)$ is a free $\CC$-module.
Denote by $\tr \colon A^{E,B,c}_n(T) \to \CC$ the internal trace function associated with the embeddings 
$A^{E,B,c}_n(T) \hra M_{\La}(\CC)$ from $\CC$-bases of $A^{E,B}_n(T)$, where
\[
\La:= L_1^2 \ldots L_n ^2.
\]
When $A^{E,B,c}_n(T)$ satisfies the assumptions of Theorem A (ii) and $L_j = d_j$, this reduces 
to the trace map in \S \ref{4.2} in the specialization $c=1$.

As in the case $c=1$, one verifies that 
 \[
z_j := c + (\ep_1 -1) y_1 x_1 + \cdots + (\ep_j -1) y_j x_j = [x_j,y_j] 
\]
are normal elements of $A^{E,B,c}_n(T)$ satisfying \eqref{z-normal}. Set $z_0=c$. Then 
\begin{equation}
\label{z-ide}
z_j = (\ep_j -1) y_j x_j + z_{j-1} \quad \mbox{and} \quad 
z_j^{d_j} = - (1 - \ep_j)^{d_j} y_j^{d_j} x_j^{d_j} +z_{j-1}^{d_j} 
\end{equation}
for $j \in [1,n]$, where the last identity is checked analogously to the one in \prref{Poisson}.
Note that $x_j^{d_j}$ and $y_j^{d_j}$ commute.

Denote $E\spcheck:=(\ep_2, \ldots, \ep_n)$. Let $B\spcheck$ be the $(n-1)\times(n-1)$ matrix obtained from $B$ by deleting 
the first row and column.
 
It follows from \leref{cent} that, if $x_1^{L_1}, y_1^{L_1}, \ldots, x_n^{L_n}, y_n^{L_n} \in \ZZ(A^{E,B,c}_n(T))$, then
\begin{equation}
\label{dL}
d_j | L_k \quad \mbox{for} \; j \leq k.
\end{equation}

\bth{general-discr} Let $A^{E,B,c}_n(T)$ be an arbitrary PI quantized Weyl algebra over an integral domain $T[c]$ 
satisfying \eqref{assum}. For a choice of central elements
\[
x_1^{L_1}, y_1^{L_1}, \ldots, x_n^{L_n}, y_n^{L_n} \in \ZZ(A^{E,B}_n(T)),
\]
denote $\AA_n:= A^{E,B,c}_n(T)$, $\CC_n := T[c, x_1^{L_1}, y_1^{L_1}, \ldots, x_n^{L_n}, y_n^{L_n}]$ and
$\AA_{n-1}:= A^{E\spcheck,B\spcheck,c'}_{n-1}(T)$, $\CC_{n-1}:= T[c', x_2^{L_2}, y_2^{L_2}, \ldots, x_n^{L_n}, y_n^{L_n}]$ for
$n>1$, $\AA_0 = \CC_0= T[c']$ for $n=1$.

{\em{(i)}} The discriminant $d(\AA_n/\CC_n)$ is a polynomial in $c^{\gcd(L_1, \ldots, L_n)}$.

{\em{(ii)}} By part {\em{(i)}} and \eqref{dL} the discriminant $d(\AA_{n-1}/\CC_{n-1})$ is a polynomial in $(c')^{d_1}$, which will be denoted by 
$d(\AA_{n-1}/\CC_{n-1})((c')^{d_1})$. We have,
\begin{align*}
d(\AA_n/\CC_n) &=_{T^\times} \theta x _1^{(L_1-d_1) \La} y_1^{(L_1-d_1) \La} (c^{L_1} - (1-\ep_1)^{L_1} y_1^{L_1} x_1^{L_1})^{(d_1-1)\La/L_1}
\\
&\times 
\prod_{i=0}^{L_1/d_1-1} \big[ d(\AA_{n-1}/\CC_{n-1})(c^{d_1} - \zeta^i (1-\ep_1)^{d_1} y_1^{d_1} x_1^{d_1})\big]^{d_1 L_1},
\end{align*}
where $\La= L_1^2 \ldots L_n^2$, $\theta = L_1^{\La} ( L_1 (1-\ep_1)^{- d_1 +1})^\La$, and 
$\zeta$ is a primitive $L_1/d_1$-st root of unity.
\eth
Note that in the setting of the theorem $\CC_n^\times = \CC_{n-1}^\times = T[c]^\times = T^\times$. Recall that by \eqref{ep1}, 
$d_1 (1-\ep_1)^{-d_1 +1} \in T$.

\bex{n12} (i) Let $n=1$. The first quantized Weyl algebra $A^{\ep_1, c}_1(T)$ is defined from $E=(\ep_1)$ and $B = (1)$ (recall that $B$ is 
a multiplicatively skewsymmetric $n \times n$ matrix). The discriminant formula in this case is
\[
d(A^{\ep_1, c}_1(T)/T[c,x_1^{L_1}, y_1^{L_1}]) =_{T^\times} \theta x_1^{(L_1-d_1) L_1^2} y_1^{(L_1-d_1) L_1^2} (c^{L_1} - (1-\ep_1)^{L_1} y_1^{L_1} x_1^{L_1})^{(d_1-1)L_1}
\]
where $\theta = L_1^{L_1^2} (L_1 (1-\ep_1)^{-d_1 +1})^{L_1^2}$.

(ii) For $n=2$, the discriminant formula is
\begin{align*}
&d(A^{E,B,c}_2(T) /T[c,x_1^{L_1}, y_1^{L_1}, x_2^{L_2}, y_2^{L_2}]) 
\\
&=_{T^\times} \theta x_1^{(L_1-d_1) \La} y_1^{(L_1-d_1) \La}  x_2^{(L_2-d_2) \La} y_2^{(L_2-d_2) \La}
\big(c^{L_1} - (1-\ep_1)^{L_1} y_1^{L_1} x_1^{L_1}\big)^{(d_1-1)\La/L_1}
\\
&\times 
\prod_{i=0}^{L_1/d_1-1} \Big[ \big(c^{d_1}- \zeta^i (1-\ep_1)^{d_1} y_1^{d_1} x_1^{d_1}\big)^{L_2/d_1} 
- (1-\ep_2)^{L_2} y_2^{L_2} x_2^{L_2} \Big]^{(d_2-1)d_1 L_1 L_2},
\end{align*}
where $\La = L_1^2 L_2^2$, $\theta= \La^{\La/2} \prod_{i=1}^2 (L_i (1- \ep_i)^{ - d_i+1})^\La$, and
$\zeta$ is a primitive $(L_1/d_1)$-st root of 1. Note that the last product in the expression for the discriminant is a polynomial in $c^{\gcd(L_1,L_2)}$.
\eex
%%%%%%%%%%%%
\subsection{Proof of \thref{general-discr}}
\label{A.2}
For a field extension $\KK'/\KK$, we will denote by 
\[
\tr_{\KK'/\KK}, \NN_{\KK'/\KK} \colon \KK' \to \KK
\]
the standard trace and norm functions. Let $\KK(\al)/\KK$ be a finite separable extension and 
$f(t)\in \KK[t]$ be the minimal polynomial of $\al$ over $\KK$. Denote by $\mu(\al)$ the $\KK$-linear endomorphism 
of $\KK(\al)$ given by multiplication by $\al$. If $\al_1=\al, \al_2, \ldots, \al_k$ are the roots of $f(t)$ 
in its splitting field, then 
\[
\mbox{the characteristic polynomial of $\mu(\al)$ is} \; \; (t- \al_1) \ldots (t-\al_k) \in \KK[t],
\]
see e.g. \cite[p. 67, Ex. 14]{R}. In particular,
\begin{equation}
\label{splitting2}
\tr \mu(\al)^j = \sum_{i=1}^k \al_i^j, \quad \NN_{\KK(\al)/\KK} (g(\al)) = \prod_{i=1}^k g(\al_i) \; \; \forall g(t) \in \KK[t].  
\end{equation}

Set
\begin{align}
\label{Delta}
\Delta &:=_{T^\times} \theta x _1^{(L_1-d_1) \La} y_1^{(L_1-d_1) \La} (c^{L_1} - (1-\ep_1)^{L_1} y_1^{L_1} x_1^{L_1})^{(d_1-1)\La/L_1}
\\
&\times 
\prod_{i=0}^{L_1/d_1-1} \big[ d(\AA_{n-1}/\CC_{n-1})(c^{d_1} - \zeta^i (1-\ep_1)^{d_1} y_1^{d_1} x_1^{d_1})\big]^{d_1 L_1},
\nn
\end{align}
\medskip
\\
\noindent
{\em{Proof of \thref{general-discr} (ii).}}
Instead of using full quantum clusters as in Sect. \ref{q-cl}, we use a part of a cluster consisting of the cluster variables 
$x_1$ and $z_1$. We localize by $x_1$ and work inductively relating the discriminant of $\AA_n$ to that of $\AA_{n-1}$.

We start with a reduction of the statement of \thref{general-discr} (ii) that has to do with the localization in question. 
This part of the theorem will follow if we establish that
\begin{equation}
\label{locDiscr}
d(\AA_n[x_1^{-L_1}]/\CC_n[x_1^{-L_1}]) =_{T[x_1^{\pm L_1}]^\times} \Delta.
\end{equation}
Indeed, if this holds then
\begin{equation}
\label{xk}
d(\AA_n/\CC_n) =_{T^\times} x_1^{k L_1} \Delta
\end{equation}
for some $k \in \Zset$. Recall the filtration from \S \ref{1.1}. By \cite[Proposition 4.10]{CPWZ2}, 
\[
\gr d(\AA_n/\CC_n)=_{T^\times} d(\gr \AA_n/ \gr \CC_n),
\]
where the second discriminant is computed with respect to the trace on $\gr \AA_n$ coming from its freeness over $\gr \CC_n$. 
Because $\gr A^{E,B,c}_n(T)$ is a localization of a skew polynomial algebra by a power of one of its generators $\ol{x}_1$, 
we can apply \cite[Proposition 2.8]{CPWZ2} to obtain
\begin{equation}
\label{grgr}
d(\gr \AA_n/\gr \CC_n) =_{T^\times} \La^\La (\ol{x}_1 \ol{y}_1)^{\La (L_1-1)} \ldots (\ol{x}_n \ol{y}_n)^{\La (L_n-1)} =_{T^\times}
\gr \De.
\end{equation}
Therefore, in \eqref{xk}, $k=0$, and \eqref{locDiscr} implies \thref{general-discr} (ii). 

Arguing as in \S \ref{4.2a}, we see that it is sufficient to prove the theorem in the case when $\ep_j - 1 \in T^\times$ for all $j \in [1,n]$. 

{\em{In the rest of the 
proof we will assume that}}
\begin{equation}
\label{ep-assume}
\ep_1 - 1 \in T^\times
\end{equation} 
{\em{and we will prove \eqref{locDiscr}.}} 
Taking associated graded in \eqref{xk} and using \eqref{grgr} gives
\[
\ol{x}^{k L_1} \gr \Delta =_{T^\times} \gr (\AA_n/\CC_n) =_{T^\times} \gr \Delta.
\]
Thus $k=0$. So, \eqref{locDiscr} implies the statement in \thref{general-discr} (ii).

The element $z_1$ commutes with all $x_j$, $y_j$ for $j>1$. The $T$-subalgebra of $\AA_n$ generated by 
$z_1$ and $x_j$, $y_j$ for $j>1$ has a natural structure of $T[z_1]$-algebra and, as such, is isomorphic 
to $\AA_{n-1}$ with $c'=z_1$. By abuse of notation we will denote this algebra by $\AA_{n-1}$ and 
its central subalgebra $T[z_1 , x_j^{L_j}, y_j^{L_j}, 2 \leq j \leq n]$ by $\CC_{n-1}$. 

Define the sets
\[
\B'':= \{ x_2^{l'_2} y_2^{l_2} \ldots x_n^{l'_n} y_n^{l_n} \mid l_j \in [0, L_j-1]\}
\]
and
\[
\B':=\{1, z_1, \ldots, z_1^{L_1 -1} \} \B'', \quad \B:= \{ 1, x_1, \ldots, x_1^{L_1 -1}\} \B', \; \B^{\opp}:= \B' \{1, x_1, \ldots, x_1^{L_1 -1} \}.
\]
By \eqref{z-ide}, 
\[
y_1^{L_1}  = \frac{(c^{d_1} - z_1^{d_1})^{L_1/d_1}}{(1-\ep_1)^{L_1} x_1^{L_1}} \cdot
\]
Therefore,
\[
\CC_n[x_1^{-L_1}] = T[c, x_1^{\pm L_1}, (c^{d_1} - z_1^{d_1})^{L_1/d_1}, x_j^{L_j}, y_j^{L_j}, 2 \leq j \leq n].
\]
Using the fact that $x_1$ normalizes $z_1$ and $x_j$, $y_j$ for $j>1$, and that $z_1$, $x_j$, $y_j$, $j>1$ generate the $T[z_1]$-algebra 
$\AA_{n-1}$, one proves that $\B$ is a $\CC_n[x_1^{-L_1}]$-basis of $\AA_n[x_1^{-L_1}]$ and $\B'$ is a $\wt{\CC}_{n-1}$-basis 
of $\AA_{n-1}[c]$, where
\[
\wt{\CC}_{n-1} = T[c, (c^{d_1} - z_1^{d_1})^{L_1/d_1}, x_j^{L_j}, y_j^{L_j}, 2 \leq j \leq n].
\]
Denote by $\tr' \colon \AA_{n-1}[c] \to \wt{\CC}_{n-1}$ the $T[c]$-linear trace function from the 
the latter basis.

We prove \eqref{locDiscr}, that is \thref{general-discr} (ii), in two steps.
\\
\noindent
{\em{Step}} I. First, we relate $d(\AA_n[x_1^{-L_1}]/\CC_n[x_1^{-L_1}])$ to $d(\AA_{n-1}[c]/\wt{\CC}_{n-1})$.
For all $b'_1, b'_2 \in \B'$ and $i,k \in [0,L_1-1]$,
\[
\tr(b'_1 x_1^i \cdot x_1^k b'_2) = \tr(b'_2 b'_1 x_1^{i+k}) = L_1 x_1^{i+k} \tr'(b'_2 b'_1) = L_1 x_1^{i+j} \tr'(b'_1 b'_2) 
\] 
if $i+j =0$ or $L_1$ and $\tr (b_1 x_1^i \cdot x_1^k b_2) =0$ otherwise. Using the standard formula for the 
determinant of a Kronecker product of matrices yields
\begin{align*}
d(\AA_n[x_1^{-L_1}]/\CC_n[x_1^{-L_1}]) &=_{T[x_1^{-L_1}]^\times} L_1^\La x_1^{(L_1-1)\La} \det( [\tr(b_1 b_2)]_{b_1 \in \B^{\opp}, b_2 \in \B})^{L_1} 
\\
&=_{T[x_1^{L_1}]^\times} L_1^\La d(\AA_{n-1}[c]/\wt{\CC}_{n-1})^{L_1}.
\end{align*} 
\noindent
{\em{Step}} II. Next, we relate $d(\AA_{n-1}[c]/\wt{\CC}_{n-1})$ to $d(\AA_{n-1} /\CC_{n-1})$. The set $\B''$ is a $\CC_{n-1}$-basis 
of $\AA_{n-1}$ (recall that $\AA_{n-1}$ and $\CC_{n-1}$ are viewed as $T[z_1]$-algebras). Denote by 
$\tr'' \colon \AA_{n-1} \to \CC_{n-1}$ the associated $T[z_1]$-linear trace and extend it to a map $\tr'' \colon \AA_{n-1}[c] \to \CC_{n-1}[c]$
by $c$-linearity. 

Denote by $\KK$ the fraction field of $T[x_1^{L_1} y_1^{L_1}]$ and let
\[
f(t) := (c^{d_1} - t^{d_1})^{L_1/d_1} - (1-\ep_1)^{L_1} x_1^{L_1} y_1^{L_1} \in \KK[t].
\]
The polynomial is irreducible, separable and $z_1$ is a root of it (for the irreducibility note that $x_1^{L_1} y_1^{L_1} \in \KK$ 
but $x_1^{d_1} y_1^{d_1} \notin \KK$). Consider the field extension $\KK(z_1)/\KK$. The traces $\tr'$ and $\tr''$ are related by 
\begin{equation}
\label{trtr}
\tr' = \tr_{\KK(z_1)/\KK} \circ \tr''  
\end{equation}
where $\tr_{\KK(z_1)/\KK}$ is extended to $\KK[c, z_1, x_j^{L_j}, y_j^{L_j}, 2 \leq j \leq n]$ by linearity on $c$ and  $x_j^{L_j}$, $y_j^{L_j}$, $j>1$.
In the proof of \eqref{trtr} we use that $z_1$ is in the center of $\AA_{n-1}[c]$.

Denote by $\al_1=z_1, \al_2, \ldots, \al_{L_1}$ the roots of $f(t)$ in its splitting field. They are given by 
\begin{equation}
\label{al-roots}
\xi^k \big( c^{d_1} - \zeta^i (1-\ep_1)^{d_1} y_1^{d_1} x_1^{d_1} )^{1/d_1} \quad \mbox{for} \; \; 
k \in [0, d_1-1], i \in [0,L_1/d_1-1],
\end{equation}
where $\zeta$ is a primitive $(L_1/d_1)$-st root of unity, as in the statement of the theorem, and
$\xi$ is a primitive $d_1$-st root of unity.

For an element $a \in \AA_{n-1}$ denote by $\tr''(a)(z_1)$ the polynomial dependence of $\tr''(a)$ on $z_1$. Combining \eqref{splitting2}
and \eqref{trtr} gives
\[
\tr'(a) = \sum_{j=1}^{L_1} \tr''(a)(\al_j). 
\]
Therefore, from the basis $\B$ and $\B'$ we have
\[
d(\AA_{n-1}[c]/\wt{\CC}_{n-1}) = \det [\tr'(z_1^{i-1} b''_1 \cdot z_1^{k-1} b''_2)]_{i,k; b''_1, b''_2} =
\det \Big[\sum_j \al_j^{i+k-2} \tr'(b''_1 b''_2)(\al_j)\Big]_{i,k; b''_1, b''_2}
\]
where in all matrices $i, j , k \in [1,L_1]$, $b''_1, b''_2 \in \B''$. The last matrix is factored into the product of block matrices
whose blocks are square matrices of size $\La/L_1^2 = | \B''|$ as follows:
\begin{equation}
\label{matrices}
[\al_j^{i-1} I ]_{i,j}  \cdot [\al_j^{i-1} I ]_{i,j} \cdot \diag( Q(\al_1), \ldots, Q(\al_{L_1}))
\end{equation}
where $i, j \in [1, L_1]$, $I$ denotes the identity matrix of size $\La/L_1^2$ and 
\[
Q(z_1):= [\tr''(b''_1 b''_2)(z_1)]_{b''_1, b''_2 \in \B''}. 
\]
By way of definition, $\det Q(z_1) =_{T\times} d(\AA_{n-1}/\CC_{n-1})(z_1^{d_1})$. This, the fact that the roots $\al_1, \ldots, \al_{L_1}$ are given by 
\eqref{al-roots}, and \thref{general-discr} (i) imply that 
the determinant of the third matrix in \eqref{matrices} equals the product in the second line of the definition \eqref{Delta} of 
$\Delta$.

It follows from \eqref{al-roots} that for all $k \in [1,L_1]$,
\[
f'(\al_j) = L_1 \al_j^{d_1-1} \frac{(1-\ep_1)^{L_1} x_1^{L_1} y_1^{L_1}}{(c^{d_1} - \al_j^{d_1})} = 
L_1 \zeta^{-i} \al_j^{d_1-1} (1- \ep_1)^{L_1 - d_1} x_1^{L_1 - d_1} y_1^{L_1 - d_1}
\]
for some $i \in [0,L_1/d_1-1]$. The determinant of the product of the first two matrices in \eqref{matrices} equals
\begin{align*}
&\prod_{1\leq i<j \leq n} (\al_i - \al_j)^{2 \La/L_1^2} = \pm \big[ \NN_{\KK(z_1)/\KK} (f'(z_1)) \big]^{\La/L_1^2} = 
\pm \prod_i f'(\al_i)^{\La/L_1^2} 
\\
&= \pm L_1^{\La/L_1} (1- \ep_1)^{(L_1 - d_1)\La/L_1} x_1^{(L_1 - d_1)\La/L_1} y_1^{(L_1 - d_1)\La/L_1} 
\big[ \NN_{\KK(z_1)/\KK}(z_1^{d_1-1}) \big]^{\La/L_1^2} 
\\
&= \pm L_1^{\La/L_1} (1- \ep_1)^{(L_1 - d_1)\La/L_1} x_1^{(L_1 - d_1)\La/L_1} y_1^{(L_1 - d_1)\La/L_1}
\big( c^{L_1} - (1- \ep_1)^{L_1} x_1^{L_1} y_1^{L_1} \big)^{(d_1-1)\La/L_1^2}
\end{align*}
using the standard expression for the discriminant of a finite separable field extension as a product of norms
\cite[pp. 66-67, Ex. 14]{Re}.
Inserting the determinants of the matrices in \eqref{matrices} in the expression for $d(\AA_{n-1}[c]/\wt{\CC}_{n-1})$, and then
using Step I and \eqref{ep-assume}, proves \eqref{locDiscr}. This completes the proof of \thref{general-discr} (ii).
\qed
\medskip
\\
\noindent
{\em{Proof of \thref{general-discr} (i).}} We prove the statement by induction on $n$: Assuming that $d(\AA_{n-1}/\CC_{n-1})$ is a 
polynomial in $z_1^L$ where $L = \gcd(L_2, \ldots, L_n)$, we show that $d(\AA_n/\CC_n)$ is a polynomial 
in $c^{\gcd(L_1, \ldots, L_n)}$. Let 
\[
d(\AA_{n-1}/\CC_{n-1})(z_1^L) = \prod_s (z_1^L - a_s)
\]
for some $a_s$ in the algebraic closure of the fraction field of $T[x_j^{L_j}, y_j^{L_j}, 2 \leq j \leq n]$. The product 
\[
\prod_{i=0}^{L_1/d_1-1} \big( (c^{d_1} - \zeta^i (1-\ep_1)^{d_1} y_1^{d_1} x_1^{d_1})^{L/d_1} - a_s \big)
\]
is a polynomial in $c^{\gcd(L_1,L)}$. This implies that the product in the 
second line of the formula for $d(\AA_n/\CC_n)$ in part (ii) of the theorem is a polynomial in $c^{\gcd(L_1,L)}$. Since the product 
on the first line of the formula is a polynomial in $c^{L_1}$, $d(\AA_n/\CC_n)$ is a polynomial in
\[
c^{\gcd(L_1, L)} = c^{\gcd(L_1, \ldots, L_n)},
\]
which proves the first part of the theorem.
\qed
%%%%%%%%%%%%%%%%%%%%%% References %%%%%%%%%%%%%%%%%%%%%%%%%%%%%%%%%%%%%%%

%%%%%%%%%%%%%%%%%%%%%%%%%%%%%%%%%%%%%%%%%%%%%%%%%%%%%%%%%%%%%%%%%%%%%%%%%%%%%%%
%%%%%%%%%%%%%%%%%%%%%%%%%%%%%%%%%%%%%%%%%%%%%%%%%%%%%%%%%%%%%%%%%%%%%%%%%%%%%%
\end{document}